 \def\simplelatex{\iftrue}
\documentclass[12pt]{article} 
\usepackage[francais]{babel}
\simplelatex%
\else%
\usepackage[latin1]{inputenc}
\usepackage[T1]{fontenc}
\fi%
\usepackage{amsmath}
\usepackage{amsfonts}
\usepackage{amssymb}
\usepackage{amsthm}
\simplelatex%
\let\mathscr=\mathcal
\else%
\usepackage{mathrsfs}
\fi%
\simplelatex%
\let\url=\texttt
\else%
\usepackage{url}
\fi%
\newenvironment{rmq}[1][]{\refstepcounter{prop}%
\bigskip
\noindent{\textbf{Remarque \theprop{}. #1}}}{}
\newenvironment{rmqs}[1][]{\refstepcounter{prop}%
\bigskip
\noindent{\textbf{Remarques \theprop{}. #1}}}{}
\newtheorem{prop}{Proposition}[section]
\newtheorem{thm}[prop]{Th\'eor\`eme}
\newtheorem{lem}[prop]{Lemme}
\newtheorem{cor}[prop]{Corollaire}

\newtheorem{conj}[prop]{Conjecture}
\newtheorem{sugn}[prop]{Suggestion} 
\newtheorem{sugns}[prop]{Suggestions}
\newtheorem{qn}[prop]{Question}

\newcommand{\Spec}{\mathop\mathrm{Spec}\nolimits}

\newcommand{\Br}{\mathop\mathrm{Br}\nolimits}

\newcommand{\Hom}{\mathop\mathrm{Hom}\nolimits}
\newcommand{\br}{\mathop\mathrm{Br}\nolimits}
\newcommand{\kker}{\mathop\mathrm{Ker}\nolimits}

\newcommand{\X}{{\mathcal X}}
\newcommand{\bP}{{\mathbb P}}
\newcommand{\Q}{{\mathbb Q}}
\newcommand{\C}{{\mathbb C}}
\newcommand{\F}{{\mathbb F}}
\newcommand{\Z}{{\mathbb Z}}
\newcommand{\N}{{\mathbb N}}
\newcommand{\R}{{\mathbb R}}
\newcommand{\A}{{\mathbb A}}
 \def\et{{\acute et}}

\simplelatex%
\let\bimu=\mu
\let\mathbi=\mathbf
\else%
\DeclareFontFamily{OML}{cmmib}{\skewchar\font127 }
\DeclareFontShape{OML}{cmmib}{m}{it}% 
       {<5><6><7><8><9>gen*cmmib%
        <10><10.95>cmmib10%
        <12><14.4><17.28><20.74><24.88>cmmib12%
        }{}
\DeclareSymbolFont{biletters}{OML}{cmmib}{m}{it}
\DeclareSymbolFontAlphabet{\mathbi}{biletters}
\DeclareMathSymbol{\bimu}{\mathord}{biletters}{"16}
\fi%
\simplelatex%

\fi%

\begin{document}
\title{Vari\'et\'es presque rationnelles, leurs points  rationnels et leurs d\'eg\'en\'erescences \\ Cours au CIME, Septembre 2007}
\author{J-L. Colliot-Th\'el\`ene}
\maketitle 
 
\section{Introduction}

Voici une s\'erie de r\'esultats classiques.

Toute forme quadratique en au moins trois variables 
 sur le corps fini $\F_{p}$ ($p$ premier) poss\`ede un z\'ero non trivial (Euler).
 Toute forme de degr\'e $d$ en $n>d$ variables sur $\F_{p}$ poss\`ede
 un z\'ero non trivial (Chevalley-Warning).
 
 Toute forme quadratique en au moins trois variables
 sur le corps $\C(t)$ des fonctions rationnelles 
en une variable poss\`ede un z\'ero non trivial (Max Noether).
Toute forme de degr\'e $d$ en   $n+1>d$ variables
 sur une extension finie  de  $\C(t)$ poss\`ede un z\'ero non trivial (Tsen).
 Ceci vaut encore sur le corps $\C((t))$  des s\'eries formelles en une variable
 (Lang).

 Sur un corps fini, sur un corps de fonctions d'une variable sur $\C$, sur le corps
 $\C((t))$, tout espace homog\`ene d'un groupe alg\'ebrique lin\'eaire connexe 
 a un point rationnel.
 
 Toute forme de degr\'e $d$ en $n>d$ variables sur le corps $p$-adique $\Q_{p}$
 poss\`ede un  z\'ero non trivial  sur une extension non ramifi\'ee de $\Q_{p}$
 (Lang).
 
 Toute forme de degr\'e $d$ en $n>d^2$ variables sur un corps de fonctions
 de deux variables sur $\C$ poss\`ede un   z\'ero non trivial  (Lang).
 
 Toute forme quadratique en  $n>2^2$ variables sur un corps $p$-adique
  poss\`ede un z\'ero non trivial (Hensel, Hasse).

 Toute forme cubique en $n>3^2$ variables sur un corps $p$-adique
 poss\`ede un z\'ero non trivial (Demjanov, Lewis).

 Pour $d$ donn\'e, pour presque tout premier $p$, toute forme
 poss\`ede un z\'ero non trivial (Ax-Kochen).

  Sur un corps $p$-adique, tout espace  homog\`ene principal d'un groupe
  semi-simple simplement connexe poss\`ede un point rationnel (Kneser, Bruhat-Tits).
 
\medskip

{\it  Sur un type donn\'e de corps, y a-t-il une classe naturelle de vari\'et\'es alg\'ebriques
 qui sur un tel  corps ont automatiquement un point rationnel ?}

 \medskip

Sur les corps de fonctions d'une variable sur $\C$ d'une part, sur les corps finis d'autre part,
des progr\`es d\'ecisifs ont \'et\'e accomplis dans les cinq  derni\`eres ann\'ees, et on peut dans une 
certaine mesure dire que la situation est  stabilis\'ee. La similitude apparente des r\'esultats
est trompeuse. Les r\'esultats cit\'es sur les corps finis s'\'etendent \`a une classe
beaucoup plus large de vari\'et\'es que les r\'esultats sur un corps de fonctions d'une variable.
Les techniques utilis\'ees sur un corps fini rel\`event de la cohomologie \'etale (ou, de la cohomologie $p$-adique).
Les techniques utilis\'ees sur un corps de fonctions sur les complexes rel\`event de
la cohomologie coh\'erente : th\'eorie de la d\'eformation, th\'eor\`emes d'annulation
de Kodaira et g\'en\'eralisations, programme du mod\`ele minimal.

Sur les corps de fonctions de deux variables, la recherche est extr\^emement active.

Dans ce rapport, qui ne contient pratiquement pas de d\'emonstrations,
j'ai essay\'e de pr\'esenter un instantan\'e de la situation.

Une partie importante du texte suit un  fil unifiant les travaux sur les corps de fonctions d'une variable, ceux sur les corps de fonctions de deux variables, et l'\'etude des vari\'et\'es sur les corps $p$-adiques.
 C'est l'\'etude des mod\`eles projectifs r\'eguliers au-dessus d'un anneau de valuation discr\`ete et de leur fibre sp\'eciale.

  \medskip

Certains aspects de ce texte ont fait l'objet  d'expos\'es depuis quelques ann\'ees.
Je remercie  Esnault, Gabber, Hassett, de Jong, Koll\'ar, Madore, Moret-Bailly, Starr et  Wittenberg pour diverses discussions.

J'engage les lecteurs \`a consulter le rapport r\'ecent d'O. Wittenberg \cite{WittS}.

\section{Notations,  rappels et pr\'eliminaires}
\label{rappels}

Soit $k$ un corps. On note $k_{s}$ une cl\^oture s\'eparable de $k$
 et  ${\overline k}$ une cl\^oture alg\'ebrique de $k$.
Une $k$-vari\'et\'e est {\it par d\'efinition} un $k$-sch\'ema s\'epar\'e de type fini sur $k$
(non n\'ecessairement irr\'eductible, non n\'ecessairement r\'eduit).
On note $X(k)=\Hom_{\Spec k}(\Spec k,X)$ l'ensemble des points $k$-rationnels 
d'un $k$-sch\'ema $X$.
Une $k$-vari\'et\'e est
dite int\`egre si elle est irr\'eductible et r\'eduite. On note alors $k(X)$ son corps des fonctions.
Une $k$-vari\'et\'e est dite g\'eom\'etriquement int\`egre si la ${\overline k}$-vari\'et\'e
$X \times_{k}{\overline k}$ est int\`egre. 
Une $k$-vari\'et\'e g\'eom\'etriquement  int\`egre poss\`ede un ouvert de Zariski
non vide qui est lisse sur $k$.
Si $k$ est un corps de caract\'eristique z\'ero, une $k$-vari\'et\'e int\`egre  $X$
est g\'eom\'etriquement int\`egre si et seulement si le corps $k$
est alg\'ebriquement ferm\'e dans le corps $k(X)$.

Pour la cohomologie galoisienne, et en particulier le groupe de Brauer d'un corps, le lecteur
consultera Serre (\cite{SerreCG}). En plusieurs endroits on fera libre usage de la notion de dimension cohomologique d'un corps.

En quelques endroits on fera aussi usage de certaines propri\'et\'es du groupe de Brauer d'un sch\'ema. 
Le lecteur se reportera aux expos\'es de  Grothendieck (\cite{GrBr}).

\begin{lem} (Nishimura, Lang)
\label{Nishi}
Soient $k$ un corps, $Z$ une $k$-vari\'et\'e r\'eguli\`ere connexe et  $Y$ une $k$-vari\'et\'e propre.
Si l'on a $Z(k) \neq \emptyset $ et s'il existe une $k$-application rationnelle de $Z$ vers $Y$, alors
$Y(k) \neq \emptyset$.
\end{lem}

\begin{lem} 
\label{Nishiirred}
Soient $k$ un corps, $Z/k$ une $k$-vari\'et\'e g\'eo\-m\'e\-tri\-quement int\`egre et
$Y/k$ une $k$-vari\'et\'e lisse connexe. S'il existe un $k$-morphisme $Z \to Y$
alors la $k$-vari\'et\'e $Y$ est g\'eo\-m\'e\-tri\-quement int\`egre.
\end{lem}
\begin{proof}  
La $k$-vari\'et\'e lisse $Y$ est g\'eo\-m\'e\-tri\-quement int\`egre si et seulement si
$Y_{k_{s}}$ est irr\'eductible. Supposons qu'elle ne le soit pas. On dispose alors
du $k_{s}$-morphisme $Z_{k_{s}} \to Y_{k_{s}}$. Le groupe de Galois de $k_{s}$
sur $k$ permute les composantes de $Y_{k_{s}}$. L'image de  $Z_{k_{s}}$
doit se trouver dans chaque composante de $Y_{k_{s}}$. Comme $Y_{k_{s}}$ est lisse, ces composantes ne se rencontrent pas.  Donc $Y_{k_{s}}$  n'a qu'une seule composante.
\end{proof}

\begin{rmq}
Comme l'observe Moret-Bailly, cet \'enonc\'e  est une cons\'e\-quence de deux r\'esultats
g\'en\'eraux.
Soit $Z \to Y$ un $k$-morphisme de $k$-vari\'et\'es.   Si $Z$ est g\'eom\'etriquement connexe et $Y$ connexe, alors $Y$
est g\'eom\'e\-tri\-quement connexe. Par ailleurs, si $Y$ est normal et
g\'eom\'etriquement connexe, alors $Y$ est  g\'eom\'etriquement  irr\'eductible.
\end{rmq}

\bigskip
\bigskip
\vfill\eject

{\bf Obstruction \'el\'ementaire}

Soient $k$ un corps, $k_{s}$ une cl\^oture s\'eparable de $k$, $\mathcal{G}={\rm Gal}(k_{s}/k)$ le groupe de Galois absolu.
Soit $X$ une $k$-vari\'et\'e lisse g\'eom\'etriquement int\`egre.  L'inclusion naturelle de groupes multiplicatifs
$k_{s}^{\times} \to k_{s}(X)^{\times} $ d\'efinit une suite exacte  $$ 1 \to  k_{s}^{\times} \to  k_{s}(X)^{\times} \to      k_{s}(X)^{\times}/  k_{s}^{\times} \to 1.$$
La classe $e(X)$ de l'extension de modules galoisiens discrets obtenue est appel\'ee l'obstruction \'el\'ementaire \`a l'existence d'un $k$-point :
si $X$ poss\`ede un $k$-point, alors $e(X)=0$
(CT-Sansuc, voir \cite{BoCTSk}).  Si $e(X)=0$, alors pour toute extension
finie s\'eparable $K/k$, l'application naturelle de groupes de Brauer $\Br K \to \Br K(X)$ est injective.

\bigskip

{\bf Construction de grands corps}
\label{grandscorps}

\medskip

Soit $k$ un corps de caract\'eristique z\'ero. Pour chaque corps $K$ contenant $k$, donnons-nous une classe ${\mathcal C}_{K}$
de $K$-vari\'et\'es alg\'ebriques g\'eo\-m\'e\-tri\-quement  int\`egres  admettant un ensemble $E_{K}$ de $K$-vari\'et\'es repr\'esentant toutes les classes de $K$-isomorphie de la classe.
Pour $k \subset K \subset L$  on suppose que le changement de corps de base
$K \to L$ envoie  ${\mathcal C}_{K}$ dans  ${\mathcal C}_{L}$.

Pour tout corps $K$ avec $k \subset K$ 
supposons satisfaite   la condition suivante :

(Stab) Si $f : X \to Y$ est un $K$-morphisme dominant
de $K$-vari\'et\'es g\'eo\-m\'e\-tri\-quement  int\`egres, si $Y$ appartient \`a ${\mathcal C}_{K}$ et si la fibre g\'en\'erique
de $f$ appartient \`a ${\mathcal C}_{K(Y)}$, alors $X$ appartient \`a ${\mathcal C}_{K}$.

Une construction bien connue, utilis\'ee  
  par Merkur'ev et Suslin (cf. \cite{Ducros}) 
   permet alors de construire un
  plongement de corps $k \subset L$ poss\'edant les propri\'et\'es suivantes :
  
  (i) Le corps $k$ est alg\'ebriquement ferm\'e dans $L$.
  
  (ii) Le corps $L$ est union de corps de fonctions de $k$-vari\'et\'es dans ${\mathcal C}_{k}$.
  
  (iii) Toute  vari\'et\'e dans ${\mathcal C}_{L}$ poss\`ede un point $L$-rationnel.
  
Le principe est le suivant : s'il existe un $k$-vari\'et\'e $X$ dans ${\mathcal C}_{k}$
qui ne poss\`ede pas de point rationnel, on remplace $k$ par le corps des fonctions
de cette vari\'et\'e. Et on it\`ere.
  Je renvoie \`a l'article de Ducros \cite{Ducros} pour la construction pr\`ecise,
  qui est reprise dans \cite{CTMadore} et \cite{CTnoteC2}.

 \medskip

Prenons pour ${\mathcal C}_{K}$ la classe des $K$-vari\'et\'es g\'eo\-m\'e\-tri\-quement int\`egres.
Rappelons qu'un corps $L$ est dit {\it pseudo-alg\'ebriquement clos} (PAC) si toute $L$-vari\'et\'e
g\'eo\-m\'e\-tri\-quement int\`egre sur $L$ poss\`ede un $L$-point.
La construction ci-dessus montre que tout corps $k$ de caract\'eristique z\'ero
est alg\'ebriquement ferm\'e dans un corps pseudo-alg\'ebriquement clos.

\medskip

En prenant pour ${\mathcal C}_{K}$ la classe des $K$-vari\'et\'es birationnelles 
\`a des fibrations successives de restrictions \`a la Weil de vari\'et\'es de Severi-Brauer,
Ducros \cite{Ducros} montre que tout corps $k$ de caract\'eristique z\'ero
est alg\'ebriquement ferm\'e dans
un corps $L$ de dimension cohomologique  $cd(L) \leq 1$.

\section{Sch\'emas au-dessus d'un anneau de valuation discr\`ete}
\label{Aschemas}

\subsection{$A$-sch\'emas de type (R),   croisements normaux, croisements normaux stricts}
\label{notation}

Soit $A$ un anneau de valuation discr\`ete, $K$ son corps des fractions, $F$ son corps r\'esiduel.
Soit $\pi$ une uniformisante de $A$.

Dans la suite de ce texte, on dira qu'un $A$-sch\'ema $\X$ {\it est de type} (R) s'il satisfait les conditions suivantes :

(i) Le  $A$-sch\'ema $\X$ est   propre  et   plat sur $A$.

(ii) Le  sch\'ema $\X$ est connexe et  r\'egulier.

(iii) La fibre g\'en\'erique $X=\X\times_{A}K=\X_{K}$ est une  $K$-vari\'et\'e
g\'eo\-m\'e\-tri\-quement int\`egre lisse sur $F$.

On note $K(X)$ le corps des fonctions de $X$,
qui est aussi celui du sch\'ema ${\X}$.
 On note $Y=\X\times_{A}F=\X_{F}$ la fibre sp\'eciale de $\X/A$.
 La fibre sp\'eciale $Y$ est le $F$-sch\'ema associ\'e  au  diviseur de Cartier
de ${\X}$ d\'efini par l'annulation de $\pi$.

Comme $\X$ est  r\'egulier donc normal,   on a une d\'ecomposition de diviseurs de Weil 
$$Y = \sum_{i} n_{i} Y_{i}$$
o\`u les $Y_{i}$ sont les adh\'erences des points $x_{i}$ de codimension 1 de $\X$ situ\'es sur la fibre sp\'eciale. L'anneau local de tout tel point $x_{i}$ est un anneau de valuation discr\`ete de corps des
fractions $K(X)$.
Si l'on note $v_{i}$ la valuation sur le corps $K(X)$ associ\'ee  \`a un tel $x_{i}$, alors
$n_{i}=v_{i}(\pi)$.

Comme ${\X}$ est  r\'egulier,   les $Y_{i}$ sont des diviseurs de Cartier sur ${\X}$.
Ce sont les composantes r\'eduites de la fibre sp\'eciale. Ce sont des $F$-vari\'et\'es  int\`egres
mais non n\'ecessairement g\'eo\-m\'e\-tri\-quement irr\'eductibles ni (si le corps $F$ n'est pas parfait)  n\'ecessairement g\'eo\-m\'e\-tri\-quement r\'eduites.

On dit que $Y \subset {\X}$ {\it est \`a croisements normaux} si  partout  localement pour la topologie
\'etale sur $ {\X}$ l'inclusion $Y \subset  {\X}$ est donn\'ee par 
une \'equation $\prod_{i \in I} x_{i}^{n_{i}}$, o\`u les $x_{i}$ font partie d'un syst\`eme
r\'egulier de param\`etres et les $n_{i}$ sont des entiers naturels.

On dit que  $Y \subset {\X}$ {\it est \`a croisements normaux stricts} 
 si  
la fibre $Y\subset {\X}$ est \`a croisements normaux et si de plus chaque composante r\'eduite
$Y_{i}$ de $Y$ est une $F$-vari\'et\'e (int\`egre) lisse. Une telle composante n'est pas n\'ecessairement
g\'eo\-m\'e\-tri\-quement irr\'eductible.

\medskip

On note $A^h$ le hens\'elis\'e de $A$, et l'on note $A^{sh}$ un hens\'elis\'e strict de $A$.
On note $K^h$ le corps des fractions de $A^h$ et $K^{sh}$ le corps des fractions 
de $A^{sh}$. Les inclusions $A \subset A^h \subset A^{sh}$ induisent
$F = F \subset F_{s}$ sur les corps r\'esiduels, o\`u $F_{s}$ est une cl\^oture s\'eparable
de $F$.

On note $\hat{A}$ le compl\'et\'e de $A$. Si les corps  $K$ et $F$ ont m\^eme caract\'eristique,
alors il existe un corps de repr\'esentants de $F$ dans $\hat{A}$ : il existe  un isomorphisme
$\hat{A} \simeq F[[t]]$.

\subsection{Quand la fibre sp\'eciale a une composante de multiplicit\'e  1}
\label{comp1}

\begin{prop}\label{multiplicite1}
Soit $\X$ un $A$-sch\'ema   de type (R).  
Les propri\'et\'es suivantes sont \'equivalentes :
\begin{itemize}
\item[(1)] 
Il existe une composante r\'eduite $Y_{i}$ dont l'ouvert de lissit\'e est non vide
et qui satisfait    $n_{i}=1$.
 \item[(2)] 
 Il existe un ouvert $U \subset \X$  lisse et surjectif sur $\Spec A$.
\item[(3)]
$\X \to \Spec A$ est localement scind\'e pour la topologie \'etale.
\item[(4)]
$\X (A^{sh}) \neq \emptyset$.
 \item[(5)] 
$X(K^{sh})  \neq \emptyset$.
\item[(6)]
$\X (\hat{A^{sh}}) \neq \emptyset$.
 \item[(7)] 
$X(\hat {K^{sh}})  \neq \emptyset$.

\end{itemize}
\label{compmult1}
\end{prop}
\begin{proof} Laiss\'ee au lecteur.
\end{proof}

Dans la situation ci-dessus, on dira que $Y$ {\it a une composante de multiplicit\'e 1}.
\footnote{La terminologie adopt\'ee dans ce texte diff\`ere de celle de \cite{BLR}. }
\begin{prop}
(a) Soient $\X$ un $A$-sch\'ema lisse connexe fid\`element plat sur $A$  
et $\X'/A$ un $A$-sch\'ema de type (R).
S'il existe une $K$-application rationnelle de $X=\X_{K}$ dans $X'=\X'_{K}$, alors
 la fibre sp\'eciale $Y'$  de $\X'/A$ a une composante de multiplicit\'e  1.
 
 (b) Soient $\X$ et $\X'$ deux $A$-sch\'emas de type (R).
  Si les fibres g\'en\'eriques $X=\X_{K}$ et $X'=\X'_{K}$ sont $K$-birationnellement
\'equivalentes, alors la fibre sp\'eciale $Y$ de $\X$ a une composante de multiplicit\'e  1 si et seulement si 
la fibre sp\'eciale $Y'$ de $\X'$  a une composante de multiplicit\'e  1.
\end{prop}
\begin{proof} Il suffit d'\'etablir (a). 
L'hypoth\`ese sur $\X/A$ et le lemme de Hensel assurent $\X(A^{sh}) \neq \emptyset$,
donc $X(K^{sh}) \neq \emptyset$.
Comme la $K$-vari\'et\'e $X$ est r\'eguli\`ere et la $K$-vari\'et\'e $X'$ propre,
d'apr\`es le lemme \ref{Nishi}
l'existence d'un $K^{sh}$-point sur $X$ implique l'existence d'un $K^{sh}$-point sur $X'$.
\end{proof}

\begin{rmq} (Wittenberg)
\label{wittr}
Soient $K=\C(u,v)$ le corps des fractions rationnelles \`a deux variables  et $X\subset \bP^2_{K}$ la conique lisse d\'efinie par l'\'equation homog\`ene
$ux^2+vy^2=z^2$. Pour tout anneau $A \subset K$ de valuation discr\`ete de rang 1, 
de corps des fractions $K$, on a $X(K^{sh}) \neq \emptyset$, o\`u
$K^{sh}$ est le corps des fractions d'un hens\'elis\'e strict $A^{sh}$ de $A$. 
Mais si  $p : \X \to S$   est un morphisme propre et plat  de vari\'et\'es
projectives lisses connexes de fibre g\'en\'erique $X/K$, le morphisme
$p$ n'est pas localement scind\'e pour la topologie \'etale.
\end{rmq}

\subsection{Quand la fibre sp\'eciale contient une sous-vari\'et\'e g\'eom\'etriquement int\`egre}
\label{compgeom}

\begin{prop}
(a) Soit $\X$ un $A$-sch\'ema r\'egulier connexe  fid\`element plat sur $A$ et soit $\X'$ un $A$-sch\'ema propre.
Si la fibre sp\'eciale $Y/F$ de $\X/A$ contient une sous-$F$-vari\'et\'e g\'eom\'etriquement int\`egre,
et s'il existe une $K$-application rationnelle de $X=\X \times_{A}K$ dans $X'=\X'\times_{A}K$,
alors la fibre sp\'eciale $Y'$ de  $\X'/A$ contient une sous-$F$-vari\'et\'e g\'eom\'etriquement int\`egre.

(b) Soient $\X$ et $\X'$ deux $A$-sch\'emas de type (R).
  Si les fibres g\'en\'eriques $X=\X_{K}$ et $X'=\X'_{K}$ sont $K$-birationnellement
\'equivalentes, alors la fibre sp\'eciale $Y$ contient une sous-$F$-vari\'et\'e g\'eom\'etriquement int\`egre si et seulement si 
la fibre sp\'eciale $Y'$ contient une sous-$F$-vari\'et\'e g\'eom\'etriquement int\`egre.
\end{prop}
\begin{proof}
Il suffit de d\'emontrer le point (a). On peut supposer la sous-vari\'et\'e int\`egre $Z \subset Y$ ferm\'ee.
Soit $Z_{lisse} \subset Z$ l'ouvert de lissit\'e de $Z/F$. Soit $p : \X_{1} \to \X$
l'\'eclat\'e de $\X$ le long de $Z$. L'image r\'eciproque de $Z_{lisse} $
dans $ \X_{1}$ est un fibr\'e projectif sur $Z_{lisse} $, qui est une $F$-vari\'et\'e g\'eom\'etriquement int\`egre. Soit $x$ son point g\'en\'erique. C'est
  un point r\'egulier de codimension 1 sur $\X_{1}$. L'application rationnelle
$\X_{1}  \to \X'$ est donc d\'efinie au point $x$. Soit $x' \in \X'$ son image.
L'adh\'erence de $x'$ dans $\X'$ est une sous-$F$-vari\'et\'e ferm\'ee de $\X'$,
munie d'une application $F$-rationnelle dominante d'une $F$-vari\'et\'e g\'eom\'etriquement int\`egre.
C'est donc une $F$-vari\'et\'e g\'eom\'etriquement int\`egre.
\end{proof}

\begin{prop}
Soit $A$ un anneau de valuation discr\`ete de corps r\'esiduel $F$ et soit $\X$ un $A$-sch\'ema   de type (R). 
Supposons les composantes r\'eduites $Y_{i}$ lisses sur $F$.
Les propri\'et\'es suivantes sont \'equivalentes :
\begin{itemize}

\item[(1)]  La fibre sp\'eciale $Y$ contient une sous-$F$-vari\'et\'e g\'eom\'etri\-quement in\-t\`egre.
\item[(2)] 
Il existe une $F$-vari\'et\'e g\'eo\-m\'e\-tri\-quement int\`egre $Z$
et un $F$-morphisme $Z \to Y$.
 \item[(3)] 
Il existe une composante r\'eduite $Y_{i}$ de $Y$ qui est g\'eo\-m\'e\-tri\-quement int\`egre.
\end{itemize}
Si de plus  $car(F)=0$, ces propri\'et\'es sont \'equivalentes aux propri\'et\'es suivantes :
\begin{itemize}
\item[(4)] Il existe une extension locale d'anneaux de valuation discr\`ete $A \subset B$
telle que le corps r\'esiduel $F=F_{A}$ de $A$ soit alg\'ebriquement ferm\'e dans le corps r\'esiduel
$F_{B}$ de $B$, et que l'on ait $X(K(B))=\X(B) \neq \emptyset$.
 \item[(5)] 
 Si $F \hookrightarrow E$ est un plongement de $F$ dans un corps pseudo-alg\'ebriquement clos $E$
   dans lequel $F$ est alg\'ebriquement ferm\'e, et si $\hat{A}=F[[t]]$,   il existe une extension
finie    totalement ramifi\'ee $L/E((t))$ avec $X(L) \neq \emptyset$.
\end{itemize}
 \label{compgeomint}
\end{prop}
\begin{proof}  Soit $f : Z \to Y$ comme en (2). Une telle application se factorise
par au moins un morphisme $Z \to Y_{i}$ pour $i$ convenable, et le lemme \ref{Nishiirred}
montre que $Y_{i}$ est alors g\'eo\-m\'e\-tri\-quement int\`egre, on a donc (3). Les autres implications entre (1), (2) et (3) sont \'evidentes.
L'\'enonc\'e (5) implique trivialement (4), et (4) implique (5) comme l'on voit en passant aux compl\'et\'es et en rempla\c cant $\hat{B} \simeq F_{B}[[u]]$ dans $E[[u]]$, o\`u $F_{B} \hookrightarrow E$
est un plongement du corps r\'esiduel $F_{B}$ dans un corps pseudo-alg\'ebriquement clos $E$
dans lequel $F_{B}$ et donc aussi $F_{A}$ est alg\'ebriquement ferm\'e.
De (4) on d\'eduit l'existence d'un $F$-morphisme $\Spec F_{B} \to Y$, 
ce qui implique l'\'enonc\'e (2).
Soit $Y_{i} \subset Y$ une composante comme en (3). 
Soit $B$ l'anneau local du point g\'en\'erique de $Y_{i} $ sur $\X$.
L'inclusion $A \subset B$ satisfait (4).
\end{proof}

Dans la situation ci-dessus, on dira que  la fibre sp\'eciale $Y/F$ {\it a une composante (r\'eduite) g\'eom\'e\-tri\-quement int\`egre}.

\bigskip

Des deux propositions pr\'ec\'edentes il r\'esulte :
\begin{prop}
 Soient $\X$ et $\X'$ deux $A$-sch\'emas de type(R), de fibres sp\'eciales respectives $Y$ et $Y'$.
 Supposons 
    les composantes r\'eduites de $Y$ et $Y'$   lisses sur $F$.
S'il existe une application $K$-rationnelle
de  $X=\X_{K}$ vers $X'=\X'_{K}$, et si 
  $Y$   a une composante  
 g\'eom\'e\-triquement int\`egre, alors
 $Y'$  a une composante g\'eo\-m\'e\-tri\-quement int\`egre.
\end{prop}

\begin{rmq}
On trouvera dans l'article  \cite{Ducros2} de Ducros  de nombreux compl\'ements et extensions des \'enonc\'es ci-dessus.
\end{rmq}

\subsection{Quand la fibre sp\'eciale a une composante
 g\'eom\'e\-triquement int\`egre de multiplicit\'e  1}
\label{compgeom1}

\begin{prop}
Soit $A$ un anneau de valuation discr\`ete de corps r\'esiduel $F$ et soit $\X$ un $A$-sch\'ema   de type (R).
 Les propri\'et\'es suivantes sont \'equivalentes :
\begin{itemize}
\item[(1)] 
Il existe une composante r\'eduite $Y_{i}$ qui est g\'eo\-m\'e\-tri\-quement int\`egre et pour laquelle
   $n_{i}=1$.
 \item[(2)] 
 Il existe un ouvert $U \subset \X$ non vide  lisse, surjectif sur $\Spec A$ et \`a fibres
 g\'eo\-m\'e\-tri\-quement int\`egres.
 \end{itemize}
 Pour $F$ de caract\'eristique z\'ero, ces propri\'et\'es sont \'equivalentes aux propri\'et\'es suivantes :
 \begin{itemize}
  \item[(3)]
  Il existe une extension non ramifi\'ee d'anneaux de valuation discr\`ete $A \subset B$
  telle que $F$ soit alg\'ebriquement ferm\'e dans le corps r\'esiduel de $B$
  et que $X(K(B))=\X(B)\neq \emptyset$.
  \item[(4)]
  Si $F \hookrightarrow E$ est un plongement de $F$ dans un corps pseudo-alg\'ebri\-que\-ment clos $E$
   dans lequel $F$ est alg\'ebriquement ferm\'e, et si $\hat{A}=F[[t]]$,
  on a $X(E((t))) \neq \emptyset$.
  \end{itemize}
L'existence d'une composante comme en (1) est une condition n\'ecessaire pour l'existence
d'un $K$-point sur $X$.
\label{compmult1geomint}
\end{prop}
\begin{proof} L'\'equivalence de (1) et (2) est claire.  L'\'equivalence de (3) et (4) est aussi claire.
Pour l'\'equivalence  entre (1) et (2) d'une part et  (3)  et (4) d'autre part, et pour la d\'emonstration de la derni\`ere assertion,  voir \cite{CTK}, fin de l'argument p.~745.
\end{proof}

Dans la situation ci-dessus, on dira que $Y$ {\it a une composante g\'eom\'e\-triquement int\`egre de multiplicit\'e  1.}

\begin{prop}
Supposons $F$ de caract\'eristique z\'ero.

(a) Soit $\X/A$ un $A$-sch\'ema connexe lisse et surjectif sur $\Spec A$,
\`a fibres g\'eom\'etriquement int\`egres. Soit $\X'/A$ un $A$-sch\'ema de type (R).
S'il existe une    application $K$-rationnelle
de  $X=\X_{K}$ vers $X'=\X'_{K}$, alors 
il existe un ouvert $U \subset \X'$ tel que le morphisme induit
$U \to \Spec A$ soit lisse et surjectif. 

(b) 
Soient $\X$ et $\X'$ deux $A$-sch\'emas de type (R).
Si les fibres g\'en\'eriques $X=\X_{K}$ et $X'=\X'_{K}$ sont $K$-biration\-nel\-lement
\'equivalentes, alors la fibre sp\'eciale $Y$ a une composante g\'eo\-m\'e\-tri\-quement int\`egre de multiplicit\'e  1 si et seulement si 
la fibre sp\'eciale $Y'$ a une composante g\'eo\-m\'e\-tri\-quement int\`egre de multiplicit\'e  1.
\end{prop}
\begin{proof}
 Cela r\'esulte imm\'ediatement du lemme \ref{Nishi} et
 de la caract\'e\-ri\-sation (3) dans la proposition \ref{compmult1geomint}.
\end{proof}

\begin{qn} 
\label{bonmodele}
Soit $k$ un corps de caract\'eristique z\'ero, $K$ un corps de type fini sur $k$,
$X$ une $K$-vari\'et\'e projective, lisse, g\'eom\'etriquement int\`egre.
Supposons que pour tout anneau de valuation discr\`ete de rang un $A$ contenant $k$
et de corps des fractions $K$ il existe un $A$-mod\`ele de type (R) de $X/K$
dont la fibre sp\'eciale contient une composante g\'eom\'etriquement int\`egre de
multiplicit\'e  1. Existe-t-il un $k$-morphisme    $\X \to  B$ de $k$-vari\'et\'es projectives, lisses,
g\'eom\'etriquement int\`egres satisfaisant les propri\'et\'es suivantes :

(a) le corps des fonctions $k(B)$ de $B$ est  $K$;

(b) la fibre g\'en\'erique de $\X \to B$ est $K$-isomorphe \`a $X$;

(c) il existe un ouvert $U \subset \X$ tel que le morphisme induit
$U \to B$ soit lisse surjectif (fid\`element plat) et \`a fibres g\'eom\'etriquement int\`egres.
\end{qn}

On comparera cette question avec la remarque \ref{wittr}.

\begin{rmq}
On exhibe facilement un morphisme  $\X \to  Y=\bP^2_{\Q}$ de $\Q$-vari\'et\'es projectives, lisses,
g\'eom\'etriquement int\`egres, de fibre g\'en\'erique une quadrique de dimension 2,
tel que la fibre en tout point de codimension 1  de $Y$ soit g\'eom\'etriquement int\`egre
sans que pour autant l'hypoth\`ese dans la question ci-dessus soit satisfaite (d\'esingulariser l'exemple de
la remarque \ref{pascodim1}).
 \end{rmq}

\subsection{Un exemple :   quadriques}

Discutons le cas des mod\`eles de quadriques de dimension au moins 1.
Supposons  $2 \in A^{\times}$. Soit $v$ la valuation de $A$.
Une quadrique lisse dans $\bP^n_{K}$ ($n \geq2$) peut \^etre d\'efinie par une forme quadratique
diagonale sur $K$. 

Consid\'erons   le cas des coniques.
En chassant les d\'enominateurs et en poussant les carr\'es
dans les variables, on voit que l'\'equation  d\'efinissant la quadrique dans $\bP^2_{K}$ peut s'\'ecrire
$$a_{0}T_{0}^2+  a_{1}T_{1}^2 + a_{2}T_{2}^2  =0    $$
avec $a_{0}, a_{1} \in A^{\times}$ et $v(a_{2})=0$ ou $v(a_{2})=1$.
Cette \'equation d\'efinit un mod\`ele r\'egulier $\X \subset \bP^2_{A}$,
et la fibre sp\'eciale $Y \subset \X$ est \`a croisements normaux.

Si $v(a_{2})=0$, alors la fibre sp\'eciale $Y/F$ est une conique lisse,
en particulier g\'eo\-m\'e\-tri\-quement int\`egre, et $Y \subset \X$
est \`a croisements normaux stricts.

Si $v(a_{2})=1$,   la fibre sp\'eciale $Y$ poss\`ede
un $F$-point rationnel \'evident, $P \in Y(F)$, donn\'e par $T_{0}=T_{1}=0$.

Si $v(a_{2})=1$ et si la classe de $-a_{0}.a_{1}$  dans $F$  est un carr\'e,
la fibre sp\'eciale se d\'ecompose sous la forme
$$Y=Y_{1}+Y_{2}$$ avec chaque $Y_{i} \simeq \bP^1_{F}$.
Dans ce cas $Y \subset \X$ est \`a croisements normaux stricts.

Si $v(a_{2})=1$ et si la classe de $-a_{0}.a_{1}$  dans $F$  n'est pas un carr\'e,
alors la fibre sp\'eciale $Y/F$ est int\`egre, mais se d\'ecompose sur une extension
quadratique de $F$ en deux droites conjugu\'ees se rencontrant en $P$,
donc $Y$ n'est pas lisse,
 $Y \subset \X$ n'est pas \`a croisements normaux stricts.
Si l'on \'eclate le point rationnel singulier $P$ sur $\X$, on obtient un mod\`ele
$\X'/A$ dont la fibre sp\'eciale $Y'$ se d\'ecompose sous la forme
$$Y'= Y'_{0} + 2E,$$
o\`u $E\subset \X$ est le diviseur exceptionnel introduit par l'\'eclatement.
La $F$-courbe $Y'_{0}$ est int\`egre, elle se d\'ecompose sur une extension quadratique
de $F$ en la somme de deux droites conjugu\'ees ne se rencontrant pas,
et rencontrant $E$ transversalement.
Donc $Y'_{0}$ est lisse, et $Y' \subset \X'$ est \`a croisements normaux stricts.

\medskip

En r\'esum\'e, pour toute conique lisse sur $K$, on a les propri\'et\'es suivantes.

(a) Il  existe un mod\`ele r\'egulier $\X$ avec $Y \subset \X$   \`a croisements normaux 
dont au moins une composante   a multiplicit\'e  1 et admet un ouvert non vide
lisse sur $F$, mais n'est pas n\'ecessairement
g\'eo\-m\'e\-tri\-quement int\`egre.

(b) Il existe un mod\`ele r\'egulier $\X$ avec $Y \subset \X$   \`a croisements normaux 
dont la fibre sp\'eciale contient une sous-$F$-vari\'et\'e g\'eo\-m\'e\-tri\-quement int\`egre.

(c) Il existe un mod\`ele r\'egulier $\X$ avec $Y \subset \X$   \`a croisements normaux stricts
dont une composante est g\'eo\-m\'e\-tri\-quement int\`egre 
(mais pas n\'eces\-saire\-ment de multiplicit\'e  1).
D'apr\`es les paragraphes \ref{comp1}  et \ref{compgeom}  les   propri\'et\'es
(a) et (b) valent pour tout $A$-mod\`ele de type (R) et la propri\'et\'e (c) vaut
pour tout $A$-mod\`ele  dont la fibre sp\'eciale est \`a croisements normaux stricts.
 
\bigskip

Consid\'erons le cas des quadriques de dimension au moins 3. On peut d\'efinir une telle
quadrique dans $\bP^n_{K} $ ($n \geq 4)$ par une \'equation
$$\sum_{i=0}^n a_{i}T_{i}^2=0,$$
avec $a_{i} \in A^{\times}$.  Dans $\bP^n_{A}$, cette \'equation  d\'efinit un mod\`ele int\`egre, normal et propre sur $A$.
La fibre sp\'eciale $Y/F$ est g\'eo\-m\'e\-tri\-quement int\`egre (et en particulier de multiplicit\'e  1).
D'apr\`es le paragraphe  \ref{compgeom1}  cette propri\'et\'e
vaut alors pour tout $A$-mod\`ele de type (R).

\section{Groupe de Brauer des sch\'emas au-dessus d'un anneau de valuation discr\`ete}
\label{Brauerdemodele}

Pour les d\'emonstrations des r\'esultats \'enonc\'es dans ce paragraphe,
le lecteur se reportera aux expos\'es de Grothendieck \cite{GrBr}.
 
 Dans cette section, la cohomologie employ\'ee
est la cohomologie \'etale, qui sur un corps est la cohomologie galoisienne du corps
(c'estt-\`a-dire de son groupe de Galois absolu).

Soit $A$ un anneau de valuation discr\`ete de corps des fractions $K$ et de corps
r\'esiduel $F$ parfait. On dispose alors d'une application r\'esidu
$$ \partial_{A} : \br K \to H^1(F,\Q/\Z)$$
envoyant le groupe de Brauer de $K$ 
dans le groupe des caract\`eres du groupe de Galois absolu de $F$.
Plus pr\'ecis\'ement, on a une suite exacte
$$ 0 \to \br A \to \br K \to H^1(F,\Q/\Z). $$
La fl\`eche de droite est surjective sur la torsion premi\`ere \`a la caract\'eristique de $F$.

Lorsque $A$ est hens\'elien, la fl\`eche naturelle $\br A \to \br F$
est un isomorphisme. Si de plus $F$ est de dimension cohomologique $\leq 1$,
alors $\br A=0$ et  $\br K \simeq  H^1(F,\Q/\Z)$. 

Soit $A \hookrightarrow B$ un homomorphisme local d'anneaux de valuation discr\`ete
\`a corps r\'esiduels parfaits.  Soit $K \subset L$ l'inclusion de corps de fractions correspondante. Soit $e$ l'indice de ramification de $B$ sur $A$,
c'est-\`a-dire la valuation dans $B$ de l'image d'une uniformisante de $A$.
Soit $F_{A} \hookrightarrow F_{B}$ l'inclusion des corps r\'esiduels. On a alors le diagramme
commutatif suivant :

$$\begin{array}{cccccccccccccccccccccccccccccc}
 \br K &\buildrel{{\partial_A}}\over \longrightarrow & H^1(F_A,{\Q}/{\Z}) \cr
\hskip8mm \downarrow{}{{\rm Res}_{K,L}}& & \downarrow{}{ e_{B/A}.{\rm Res}_{F_A,F_B}}\cr
\br L &\buildrel{{\partial_B}}\over \longrightarrow &H^1(F_B,{\Q}/{\Z}).\cr
\end{array} $$ 
 
Soit $F'_{B} \subset F_{B}$ la fermeture alg\'ebrique de $F_{A}$ dans $F_{B}$.
 Le noyau de $$e_{B/A}.{\rm Res}_{F_A,F_B} : H^1(F_A,\Q/\Z) \to H^1(F_B,\Q/\Z)$$
s'identifie au noyau de 
$$e_{B/A}.{\rm Res}_{F_A,F_B} : H^1(F_A,\Q/\Z) \to H^1(F'_B,\Q/\Z)$$

Soit $A$ un anneau de valuation discr\`ete de corps des fractions $K$
de corps r\'esiduel un corps $F$ de caract\'eristique z\'ero.
Soit  $\X$ un $A$-sch\'ema de type (R). Soit $X/K$ la fibre g\'en\'erique.
Soit $Y=\sum_{i}e_{i} Y_{i}$ la d\'ecomposition de la fibre sp\'eciale  $Y$ en
diviseurs int\`egres. Soit $F_{i}$ la fermeture alg\'ebrique de $F$ dans $F(Y_{i})$.

Comme les sch\'emas int\`egres $\X$ et $X$ sont r\'eguliers, les applications
de restriction $\br \X \to \br X \to \br K(X)$ sont injectives.
On dispose alors du diagramme commutatif de suites exactes
$$\begin{array}{ccccccccccccccccccccccccccccccccccccccccccccccccccc}
0 & \to & \br A  & \to &  \br K &\buildrel{{\partial_A}}\over \longrightarrow  & H^1(F,{\Q}/{\Z}) \cr
&   & \downarrow &&\downarrow{}{    }     & & \downarrow{}{    e_{i}. {\rm Res}_{F,F(  Y_{i}   )}     } \cr
0 & \to & \br \X & \to & \br X  &\buildrel{ \oplus_{i}{\partial_i}}\over \longrightarrow  & \oplus_{i}   H^1(F(Y_{i}),\Q/\Z).     \cr
\end{array} $$ 
et de la suite exacte qui s'en d\'eduit
$$ 0 \to \kkerÊ[\br A \to \br K(X)] \to \kker [\br K \to \br K(X)] \to \hfill$$
$$ \hskip5cm   
  \kker [ H^1(F,{\Q}/{\Z}) \buildrel{ \oplus_{i} e_{i}.{\rm Res}_{F,F_{i}} 
}\over \longrightarrow    \oplus_{i}   H^1(F_{i},\Q/\Z)]$$

 \begin{prop}
 \label{injbrfibres}
 Soit $A$  un anneau de valuation discr\`ete hens\'elien \`a corps r\'esiduel
 $F$ de caract\'eristique z\'ero et de dimension cohomologique au plus 1.
 Soit $\X$ un $A$-sch\'ema de type (R), de fibre g\'en\'erique $X$.
 Avec les notations ci-dessus, les deux  propri\'et\'es suivantes sont \'equivalentes :
 
 (i) L'application $\Br K \to \Br X/ \Br \X$ est injective.
 
 (ii) L'application
 $ H^1(F,{\Q}/{\Z}) \buildrel{ \oplus_{i} e_{i}.{\rm Res}_{F,F_{i}} 
}\over \longrightarrow    \oplus_{i}   H^1(F_{i},\Q/\Z) $
est injective.

En particulier, si la fibre sp\'eciale $Y$ poss\`ede une composante g\'eom\'e\-tri\-quement int\`egre
de multiplicit\'e  1, ou plus g\'en\'eralement si le pgcd des entiers  $e_{i}.[F_{i}:F]$
est \'egal \`a 1, alors  $\Br K \to \Br X/ \Br \X$   est injective, et il en est donc
de m\^eme de $\Br K \to \Br X$ et de $\Br K \to \Br K(X)$.
   \end{prop}

Je renvoie \`a  \cite{CTSaito} pour une discussion du cas o\`u le corps
r\'esiduel $F$ est fini.

 \section{Corps $C_{i}$}

\begin{thm}(Tsen, 1933) 
\label{Tsen}
Soit $F$ un corps alg\'ebriquement clos et $K=F(C)$ un corps de fonctions d'une variable
sur $F$. Soit $X \subset \bP^n_{K}$ une hypersurface de degr\'e $d$. Si l'on a $n \geq d$,
alors $X(K)\neq \emptyset$.
\end{thm}

On notera que l'on ne fait aucune hypoth\`ese sur $X$, qui peut  \^etre r\'eductible.

Le cas des coniques ($d=2, n=2$) avait \'et\'e \'etabli  par Max Noether par une m\'ethode
g\'eom\'etrique.

\medskip
Soit $i \geq 0$ un entier. 
On dit qu'un corps $K$ poss\`ede la propri\'et\'e $C_{i}$ si
toute forme \`a coefficients dans $K$, de degr\'e $d$ en $n+1 > d^i$ variables a un z\'ero non trivial
sur $K$.
On dit qu'un corps $K$ poss\`ede la propri\'et\'e $C'_{i}$ si  pour toute famille finie de formes
$\{\Phi_{j}(X_{0}, \dots,X_{n})\}_{  j=1,\dots,r}$  de degr\'es respectifs $d_{1}, \dots, d_{r}$
avec $n+1 > \sum_{j=1}^r d_{j}^{i}$ il existe un z\'ero commun non trivial sur $K$.

La propri\'et\'e $C_{i}$ implique la propri\'et\'e ci-dessus pour
un syst\`eme de formes $\{\Phi_{j}\}$  lorsque tous les degr\'es $d_{j}$
sont \'egaux (Artin, Lang, Nagata). 
On ne sait pas si en g\'en\'eral $C_{i}$ implique $C'_{i}$
(voir \cite{Pf}).

Un corps est alg\'ebriquement clos si et seulement si il est  $C_{0}$.
Le th\'eor\`eme de Tsen dit qu'un corps de fonctions d'une variable sur
un corps alg\'ebriquement clos est un corps $C_{1}$.
Dans sa th\`ese, suivant des suggestions d'E.~Artin, 
S.~Lang 
g\'en\'eralisa le th\'eor\`eme de Tsen. Son r\'esultat, pour lequel on trouve
quelques ant\'ec\'edents dans les textes des g\'eom\`etres alg\'ebristes  italiens, fut compl\'et\'e par Nagata. 
Le r\'esultat g\'en\'eral est le suivant.

\begin{thm}(Lang, Nagata) \cite{L}
\label{LangCi}
Soit $K$ un corps $C_{i}$. 
Toute extension alg\'ebrique de $K$ est un corps $C_{i}$.
Le corps des fractions rationnelles en une variable $K(t)$ est $C_{i+1}$.
De fa\c con g\'en\'erale,
toute extension de degr\'e de transcendance
$n$ de $K$ est un corps $C_{i+n}$.
\end{thm}

Lang \'etablit aussi le th\'eor\`eme suivant.

\begin{thm}(Lang)  \cite{L}
\label{LangWitt}
Soit $A$ un anneau de valuation discr\`ete hens\'elien de corps des fractions $K$
et de corps r\'esiduel $F$. Soit  $\hat{K}$ le compl\'et\'e de $K$. Supposons $\hat{K}$ s\'eparable sur $K$.
Si $F$ est alg\'ebriquement clos, alors $K$ est un corps $C_{1}$.
\end{thm}

En particulier l'extension maximale non ramifi\'ee $\Q_{p}^{nr}$ du corps
$p$-adique $\Q_{p}$ est un corps $C_{1}$.

Dans la situation consid\'er\'ee au paragraphe \ref{Aschemas} (voir la Proposition \ref{multiplicite1}),  
si l'on suppose le corps r\'esiduel $F$
de $A$ parfait, ce th\'eor\`eme assure que la fibre sp\'eciale d'une hypersurface
de degr\'e $d$ dans $\bP^n_{K}$ avec $n \geq d$ poss\`ede une composante de multiplicit\'e~1.

\begin{thm}(Greenberg)\cite{Greenberg}
\label{Greenberg}
 Si $F$ est un corps $C_{i}$  alors $K=F((t))$ est un corps
$C_{i+1}$.
\end{thm}

Ce th\'eor\`eme ne s'\'etend pas dans une situation d'in\'egale caract\'eristique : les corps
$p$-adiques ne sont pas $C_{2}$ (Terjanian). Ils ne sont en fait    $C_{i}$ pour aucun $i$
(Arkhipov et Karatsuba).

 \medskip
 
 Le th\'eor\`eme de Tsen est souvent mis en parall\`ele avec l'\'enonc\'e suivant,
 qui implique que les corps finis sont des corps $C_{1}$.
\begin{thm} (Chevalley, Warning, 1935)
\label{chevalleywarning}
Soit $\F$ un corps fini de caract\'eristique $p$. Soit $X \subset \bP^n_{\F}$ une hypersurface de degr\'e $d$. Si l'on a $n \geq d$,
alors $X(\F)\neq \emptyset$. Plus pr\'ecis\'ement, le cardinal de $X(\F)$ est congru \`a 1 modulo $p$.
\end{thm}

On montra plus tard (Ax (1964), Katz (1971)) que si le cardinal de $\F$ est $q$ alors
  le cardinal de $X(\F)$ est congru \`a 1 modulo $q$.  

Le cas des coniques avait \'et\'e \'etabli par Euler \cite{Euler}, dans un article  o\`u il \'etablit aussi  la formule
de multiplication pour les sommes de quatre carr\'es. La combinaison de ces deux r\'esultats
lui permit de montrer que tout rationnel positif est une somme de quatre carr\'es de rationnels.

Comme pour le th\'eor\`eme de Tsen, le th\'eor\`eme de Chevalley-Warning ne fait aucune hypoth\`ese sur $X$.
On peut voir l\`a l'origine de la conjecture suivante.

\begin{conj} \label{Ax} (Ax) \cite{Ax} Soient $K$ un corps et $X \subset \bP^n_{K}$ une hypersurface de degr\'e $d$.
Si l'on a $n \geq d$, alors il existe  une sous-$K$-vari\'et\'e $Y\subset X$ qui est g\'eo\-m\'e\-tri\-quement irr\'eductible.
\end{conj}

Si $K$ est parfait, dans la conclusion on peut remplacer  {\og g\'eo\-m\'e\-tri\-quement irr\'eductible\fg}    par {\og g\'eo\-m\'e\-tri\-quement int\`egre \fg}. Mais comme l'exemple de la forme irr\'eductible $T_{0}^2+xT_{1}^2+yT_{2}^2$
sur le corps $K=\F_{2}(x,y)$ le montre,  ceci ne vaut pas sur $K$ corps non parfait.

Le cas $d=2$ a \'et\'e discut\'e plus haut.
Le cas $d=3$ est facile. Le cas $d=4$ fut \'etabli par Denef, Jarden et Lewis dans \cite{DJL}.
Dans le m\^eme article, 
les auteurs \'etablissent la conjecture lorsque $K$
 contient un corps alg\'ebriquement clos. La d\'emonstration de ce r\'esultat utilise
 la th\'eorie des corps hilbertiens.
  
 En caract\'eristique nulle, la conjecture d'Ax est maintenant un th\'eor\`eme de Koll\'ar   (Th\'eor\`eme  \ref{KollarAx} ci-apr\`es).

 \medskip
 
 Les corps finis, les corps de fonctions d'une variable, le corps $\C((t))$  sont des corps de dimension $\leq 1$  au sens de Serre (\cite{SerreCG}, II.3.1)  : le groupe de Brauer de toute extension finie de $k$ est trivial. De fait, tout corps $C_{1} $ est de dimension  $\leq 1$. On sait (Ax) que la r\'eciproque est fausse (pour  des r\'ef\'erences et d'autres r\'esultats dans cette direction, voir \cite{CTMadore}).

Tout espace homog\`ene d'un groupe alg\'ebrique lin\'eaire connexe sur un corps parfait de dimension $\leq 1$ a un point rationnel (Steinberg, Springer, voir \cite{SerreCG}).

Les corps finis et le corps $\C((t))$ ont des groupes de Galois isomorphes au groupe $\hat{\Z}$.
Tout espace homog\`ene d'une vari\'et\'e ab\'elienne sur un corps fini poss\`ede un point rationnel (Lang) mais ceci n'est pas vrai sur $\C((t))$, comme le montre l'exemple de la courbe de genre 1
donn\'ee dans $\bP^2_{\C((t))}$ par l'\'equation homog\`ene $X^3+tY^3+t^2Z^3=0$.
 
\section{$R$-\'equivalence et \'equivalence rationnelle sur les z\'ero-cycles}

 \medskip
 
 Soient $k$ un corps et $X$ une $k$-vari\'et\'e. Deux $k$-points de $X$ sont \'el\'e\-men\-tai\-rement $R$-li\'es s'il existe un $k$-morphisme $U \to X$ d'un ouvert $U$ de $\bP^1_{k}$ tel que les deux points soient dans l'image de $U(k)$. La relation d'\'equivalence sur $X(k)$ engendr\'ee par cette relation est appel\'ee la  $R$-\'equivalence \cite{Manin}. Pour $k$ de caract\'eristique z\'ero,
 l'ensemble $X(k)/R$ est un invariant $k$-birationnel des $k$-vari\'et\'es projectives, lisses, g\'eo\-m\'e\-tri\-quement int\`egres \cite{RET}.

La $R$-\'equivalence a \'et\'e beaucoup \'etudi\'ee lorsque $X$ est un $k$-groupe lin\'eaire \cite{RET, GilleIHES}.  

 \medskip
   
   Soit $X/k$ une $k$-vari\'et\'e alg\'ebrique. On note $Z_{0}(X)$ le groupe ab\'elien libre sur les
   points ferm\'es de $X$. C'est le groupe des z\'ero-cycles de $X$.
   Le groupe de Chow (de degr\'e z\'ero), not\'e  $CH_{0}(X)$, est le quotient du groupe $Z_{0}(X)$
   par le sous-groupe engendr\'e par les z\'ero-cycles de la forme $\pi_{*}({\rm div}_{C}(g))$,
   o\`u $\pi :C \to X$ est un $k$-morphisme propre d'une $k$-courbe $C$ normale int\`egre,
     $g$ est une fonction rationnelle sur $C$ et ${\rm div}_{C}(g) \in Z_{0}(C)$ est son diviseur.
   
   \medskip
   
   Si $X/k$ est propre, l'application lin\'eaire   ${\rm deg}_{k} : Z_{0}(X) \to \Z$  envoyant un point ferm\'e
   $P$ sur son degr\'e $[k(P):k]$ passe au quotient par l'\'equivalence rationnelle ci-dessus d\'efinie.
   On note   $A_{0} (X)$ le noyau de l'application ${\rm deg}  : CH_{0}(X) \to \Z$, et on l'appelle
   le groupe de Chow r\'eduit de $X$.
   
   Le groupe de Chow r\'eduit  est un invariant $k$-birationnel des $k$-vari\'et\'es projectives, lisses,
   g\'eo\-m\'e\-tri\-quement int\`egres.
      
   \medskip
   
   Pour $X/k$ propre,  l'application \'evidente  $X(k) \to Z_{0}(X)$
induit une application
   $$ X(k)/R \to CH_{0}(X)$$
   dont l'image tombe dans l'ensemble des classes de cycles de degr\'e 1.

\section{Autour du th\'eor\`eme de Tsen : Vari\'et\'es rationnellement connexes}

Dans le programme de classification de Mori est apparue au d\'ebut des ann\'ees 1990
la notion de vari\'et\'e rationnellement connexe. Les travaux fondateurs r\'esultent 
d'une  collaboration entre
Koll\'ar, Miyaoka et Mori;  certains des r\'esultats sont dus  \`a Campana.
Un r\^ole-cl\'e y est jou\'e par la th\'eorie des d\'eformations, plus pr\'ecis\'ement par l'\'etude infinit\'esimale des sch\'emas $\rm{Hom}$, cas particulier des sch\'emas de Hilbert.
 On consultera  les livres \cite{Kollarlivre} et \cite{Debarre},  ainsi que les articles \cite{AraujoKollar}
 et le   r\'ecent rapport \cite{Starrrapport}.

Tout fibr\'e vectoriel sur la droite projective est isomorphe \`a une somme directe
de fibr\'es de rang 1, donc de la forme ${\mathcal O}(n), n \in \Z$.
Soit $k$ un corps alg\'ebriquement clos.
Soit $X$ une $k$-vari\'et\'e alg\'ebrique projective, lisse, connexe, de dimension $d$.
Soit $T_{X}$ son fibr\'e tangent.
On dit que $X$ est s\'eparablement rationnellement connexe (SRC) s'il existe un
morphisme $f : \bP^1 \to X$ tr\`es libre, c'est-\`a-dire tel que dans une (donc dans  toute) d\'ecomposition
du fibr\'e vectoriel $$f ^*T_{X} \simeq \oplus_{i=1}^d  {\mathcal O}(a_{i}),$$  on ait $a_{i} \geq 1$ pour tout $i$.

On dit qu'une vari\'et\'e projective lisse et connexe $X$ est rationnellement connexe (RC) si, pour tout corps alg\'ebriquement clos $\Omega$
contenant $k$, par un couple g\'en\'eral 
de points de  $X(\Omega)$ il passe une courbe de genre z\'ero, i.e. il existe un $\Omega$-morphisme
$\bP^1_{\Omega} \to X_{\Omega}$ dont l'image contient les deux points.

On dit qu'une vari\'et\'e projective lisse et connexe $X$ est rationnellement connexe par cha\^{\i}nes (RCC)  si, pour tout corps alg\'ebriquement clos $\Omega$
contenant $k$, tout couple g\'en\'eral  de points  de $X(\Omega)$  est li\'e  par une cha\^{\i}ne de courbes de genre z\'ero.
 Cette derni\`ere propri\'et\'e est \'equivalente \`a la condition que  tout couple de points de  $X(\Omega)$ est li\'e par une cha\^{\i}ne de courbes de genre z\'ero. En d'autres termes,
 l'ensemble $X(\Omega)/R$ est r\'eduit \`a un \'el\'ement.
 
 Au lieu de faire les hypoth\`eses ci-dessus pour tout corps alg\'ebriquement clos $\Omega$
 contenant $k$, il suffit de les faire pour un tel corps non d\'enombrable.
 
 Toute vari\'et\'e    RC est clairement RCC.
 
 \begin{thm} (Koll\'ar-Miyaoka-Mori)
 Toute vari\'et\'e 
   SRC est RC donc RCC. En caract\'eristique z\'ero, ces trois propri\'et\'es
 sont \'equiva\-lentes, et elles impliquent :
 
 Pour tout corps alg\'ebriquement clos $\Omega$
contenant $k$ et tout ensemble fini de points $x_{1}, \dots, x_{n} \in X(\Omega)$
il existe un morphisme $f : \bP^1_{\Omega} \to X_{\Omega}$ tr\`es libre
tels que tous les $ x_{i}$ soient dans l'image de $f$.
 \end{thm}
 
 En dimension 1, une vari\'et\'e 
  est  RC si et seulement si elle est une courbe lisse de genre z\'ero.
En dimension 2, une vari\'et\'e 
est SRC si et  seulement si elle est  rationnelle, i.e. birationnelle \`a un espace projectif.
Une vari\'et\'e  projective et lisse   unirationnelle est RC. Une vari\'et\'e projective et lisse  s\'eparablement unirationnelle est SRC.
 Sur $k$ alg\'ebriquement clos, le groupe de Chow r\'eduit $A_{0}(X)$ d'une vari\'et\'e RCC
 est clairement trivial.
 
Si $k$ est un corps quelconque, une $k$-vari\'et\'e  
est dite rationnellement connexe resp. rationnellement con\-nexe par cha\^{\i}nes, respectivement s\'eparablement rationnellement connexe  si elle est  g\'eom\'etriquement int\`egre et si elle est
RC, resp. RCC, resp. SRC,  apr\`es passage \`a un corps alg\'ebriquement clos contenant $k$.

Les compactifications lisses d'espaces homog\`enes de groupes alg\'ebriques lin\'eaires connexes sont des vari\'et\'es  RC.

En dimension 2, on dispose d'une classification $k$-birationnelle des $k$-surfaces SRC, c'est-\`a-dire
des $k$-surfaces projectives, lisses (g\'eo\-m\'e\-tri\-quement) rationnelles (Enriques, Manin, Iskovskikh, Mori) : toute telle surface est
$k$-birationnelle \`a une $k$-surface de del Pezzo ou \`a une $k$-surface fibr\'ee en coniques
au-dessus d'une conique lisse. Une surface de del Pezzo $X$ est une surface projective et lisse
dont le fibr\'e anticanonique $\omega_{X}^{-1}$ est ample. Le degr\'e d'une telle surface est
l'entier $d =(\omega.\omega)$. Il satisfait $1 \leq d \leq 9$.

Le groupe de Chow r\'eduit $A_{0}(X)$ d'une vari\'et\'e  
RCC sur un corps $k$ quelconque est clairement un groupe de torsion. On peut montrer (Prop. \ref{CTinv06}) qu'il est annul\'e par un entier $N=N(X)>0$.

Une vari\'et\'e de Fano est une vari\'et\'e projective lisse dont fibr\'e anticanonique est ample.
Si $X \subset \bP^n$ est une intersection compl\`ete lisse connexe d\'efinie par
des formes $\Phi_{j}, j=1,\dots,r$ de degr\'es respectifs $d_{j}, j=1,\dots,r$, 
alors $X$ est de Fano si et seulement si $n  \geq \sum_{j}d_{j}$.
On reconna\^{\i}t l\`a la condition $C_{1}$.

\medskip
Un th\'eor\`eme difficile est le suivant  (\cite{Camp}, \cite{KMM}; voir aussi \cite{Kollarlivre}, \cite{Debarre}) :
\begin{thm} (Campana, Koll\'ar-Miyaoka-Mori)
Une vari\'et\'e de Fano est rationnellement connexe par cha\^{i}nes.
\end{thm}

Le th\'eor\`eme suivant peut donc \^etre vu comme une g\'en\'eralisation du th\'eor\`eme de Tsen.
 
\begin{thm} (Graber, Harris, Starr \cite{GHS} ; de Jong, Starr \cite{dJS})
\label{GHSJ}
Soit $F$ un corps alg\'ebriquement clos et $K=F(C)$ un corps de fonctions
d'une variable sur $F$. Soit $X$ une $K$-vari\'et\'e  
 s\'eparablement rationnellement connexe.
Alors $X(K) \neq \emptyset$.
\end{thm}

Ce th\'eor\`eme implique le r\'esultat suivant.
\begin{thm}
\label{ratconsurratcon}
Soit $F$ un corps alg\'ebriquement clos de caract\'eristique z\'ero.
Soit $f : X \to Y$ un morphisme dominant de  $F$-vari\'et\'es projectives
et lisses,  \`a fibre g\'en\'erique g\'eom\'etriquement int\`egre.
  Si $Y$ est rationnellement connexe et si la fibre g\'en\'erique
est une vari\'et\'e rationnellement connexe, alors $X$
est une vari\'et\'e rationnellement connexe.
\end{thm}

Le th\'eor\`eme \ref{GHSJ} a aussi le corollaire suivant, connu des experts.
\begin{thm}
\label{SRCsurF((t))}
Soit $R$ un anneau de valuation discr\`ete hens\'elien
\'equi\-carac\-t\'eristique de corps r\'esiduel $F$ alg\'ebriquement clos, de corps des fractions $K$. Soit $X \subset \bP^n_K$ une $K$-vari\'et\'e 
s\'eparablement rationnellement connexe.
Alors $X(K) \neq \emptyset$.
\end{thm}

\begin{proof}  
 Soit $S $ le compl\'et\'e de $R$ et $L $ le corps des fractions de $S$.
 Comme $X/K$ est lisse, $X(K)$ est dense dans $X(L)$ pour la topologie d\'efinie
 par la valuation de $K$  (\cite{BLR},  Chap. 3.6, Cor. 10 p.~82). Il suffit donc d'\'etablir
 le th\'eor\`eme en rempla\c cant $R$ par $S$ et $K$ par $L$.
 L'anneau $S$ admet alors un corps de repr\'esentants isomorphe \`a $F$,
 on peut donc identifier $S=F[[t]]$ et  $L=F((t))$.
 
 Les lettres $R$ et $K$  \'etant d\'esormais libres, notons maintenant $R$  
  le hens\'elis\'e de $F[t]$ en $t=0$ et $K$  le corps des fractions de $R$.
  Le corps $L=F((t))$ est le compl\'et\'e de $K$.
 
Soit $\X' \subset  \bP^n_{S}$ l'adh\'erence sch\'ematique de $X \subset  \bP^n_{L}$.
C'est un sch\'ema int\`egre projectif et plat sur $S$.
 La $F$-alg\`ebre $S $ est la limite inductive filtrante de ses $F$-sous-alg\`ebres de type fini.
Il existe donc 
une $F$-alg\`ebre de type fini int\`egre $A \subset S$ et un $A$-sch\'ema $\X$ projectif et plat
  tel que $\X\times_{A}S \simeq \X'$.
En particulier la fibre g\'en\'erique de $\X/A$ est une vari\'et\'e s\'eparablement rationnellement connexe.
Il existe un  ouvert non vide $U \subset Y=\Spec A$, qu'on peut prendre lisse sur $F$,  tel que toutes les fibres du morphisme $\X \to Y$
  au-dessus de $U$  sont   des vari\'et\'es s\'eparablement rationnellement connexes 
   (Koll\'ar, Miyaoka, Mori, voir  \cite{Kollarlivre} IV.3.11).

Notons $\xi \in Y(S) \subset Y(L)   $ le point correspondant \`a $A \subset S$.
Munissons $U(L)$ de la topologie d\'efinie par la valuation de $L$.
L'ensemble $U(K)$ est dense dans $U(L)$    (\cite{BLR},  Chap. 3.6, Cor. 10 p.~82). 
Pour tout entier $n\geq 1$ l'application naturelle $Y(S) \to Y(S/t^n)$ a ses fibres ouvertes
dans $Y(S) Ê\subset Y(L)$.
On peut donc trouver dans $U(K) \subset Y(K)$ un point $\xi_{n}$ qui soit dans
$Y(S)$ et donc dans $Y(R)$ et qui ait m\^eme image que $\xi$ dans $Y(R/t^n)=
Y(S/t^n)$.
 On peut donc pour tout $n\geq 1$ trouver un $R$-sch\'ema projectif  
 $\X_{n}$
  \`a fibre g\'en\'erique $X_{n}/K$ s\'eparablement rationnellement connexe,
  tel que $\X_{n} \times_{R} R/t^n \simeq \X\times_{A} \times_{S}S/t^n$
  (noter que l'on a $R/t^n=S/t^n$.)
Le corps des fractions du hens\'elis\'e $R$ de $F[t]$ en $t=0$ est
l'union de corps de fonctions de $F$-courbes.
 Le th\'eor\`eme \ref{GHSJ}
 assure donc $X_{n}(K)\neq \emptyset$. Comme $\X_{n}/R $ est projectif,
 on a donc $\X_{n}(R) =X_{n}(K) \neq \emptyset$, donc $\X(S/t^n)=\X_{n}(R/t^n) \neq \emptyset$.
  On a donc $\X(S/t^n) \neq \emptyset$ pour tout entier $n$.
   Par  un th\'eor\`eme de
 Greenberg (\cite{Greenberg}, Thm. 1)
 ceci implique
  $\X(S) \neq \emptyset$. On a donc $X(L) \neq \emptyset$.
\end{proof}

On ne conna\^{\i}t pas de d\'emonstration de ce th\'eor\`eme qui ne passe
pas par le cas global $K=F(C)$.

 \medskip
Le th\'eor\`eme \ref{GHSJ} admet
un  th\'eor\`eme {\og r\'eciproque \fg}  : 

\medskip
\begin{thm} (Graber-Harris-Mazur-Starr)\cite{GHMS}
 Soit $k=\C$. Soit $S$ une vari\'et\'e lisse sur $\C$ de dimension au moins 2. Soit $X \to S$
un morphisme projectif et lisse \`a fibre g\'en\'erique g\'eo\-m\'e\-tri\-quement int\`egre.
Si la  restriction de $X \to S$ \`a toute courbe $C \subset S$ admet une section, 
alors il existe une $\C(S)$-vari\'et\'e $Z$ g\'eo\-m\'e\-tri\-quement int\`egre et rationnellement connexe
et un $\C(S)$-morphisme de $Z$ dans la fibre g\'en\'erique de $X \to S$.
\end{thm}

\medskip

Les vari\'et\'es rationnellement connexes sont donc en quelque sorte caract\'eris\'ees par
le fait d'avoir automatiquement un point sur le corps des fonctions d'une courbe sur les complexes.

La classe des vari\'et\'es   rationnellement connexes semble \^etre la classe la plus large
de vari\'et\'es projectives lisses
\`a laquelle on peut \'etendre le th\'eor\`eme de Tsen. 

Un exemple explicite de surface d'Enriques
sur $K=\C((t))$ sans $K$-point  a \'et\'e construit par Lafon \cite{Lafon}.
Un mod\`ele affine, avec variables $x,y,u,z,$ est d\'efini par le syst\`eme
    $$x^2-tu^2+t=(t^2u^2-t)y^2$$ $$ x^2-2tu^2+(1/t)=t(t^2u^2-t)z^2.$$
     
[Dans la classification des surfaces, les surfaces d'Enriques sont en quelque sorte les plus proches des 
surfaces rationnelles, une telle surface $X$ satisfait en particulier $H^1(X,O_{X})=0$ et $H^2(X,O_{X})=0$.]

\medskip

En caract\'eristique z\'ero,
le th\'eor\`eme \ref{SRCsurF((t))} est \'equivalent \`a l'assertion suivante : 
 \begin{thm}
\label{fibSRCm1}
Soit $A$ un anneau de valuation discr\`ete de corps r\'esiduel de caract\'eristique z\'ero
et soit $\X$ un $A$-sch\'ema int\`egre propre r\'egulier. Si la  fibre g\'en\'erique
est une vari\'et\'e SRC,  alors la fibre sp\'eciale poss\`ede une composante
de multiplicit\'e  1.
 \end{thm}
 
Comme indiqu\'e ci-dessus, le th\'eor\`eme ci-dessus ne vaut d\'ej\`a plus
si la fibre g\'en\'erique est une surface d'Enriques.

\begin{rmq}
On ne sait pas si l'analogue de ce th\'eor\`eme vaut dans le cas d'in\'egale caract\'eristique.
Par exemple, si $X$ est une vari\'et\'e rationnellement connexe sur le corps $p$-adique
$\Q_{p}$, a-t-elle un point dans une extension non ramifi\'ee de $\Q_{p}$ ?
C'est  vrai en dimension 1 ou 2. En effet l'extension maximale non ramifi\'ee de $\Q_{p}$
est un corps $C_{1}$ (th\'eor\`eme de Lang) et par inspection de la classification
$k$-birationnelle des $k$-surfaces rationnelles, on montre (Manin, l'auteur) que
toute $k$-surface rationnelle (projective et lisse) sur un corps $k$ qui est $C_{1}$ poss\`ede un point $k$-rationnel.
\end{rmq}

 Motiv\'e par la conjecture \ref{Ax} (Ax),  par les \'enonc\'es  des th\'eor\`emes
 de Tsen et de Chevalley-Warning,
 et par plusieurs r\'esultats qui seront discut\'es plus bas,
 on peut, suivant Koll\'ar \cite{KollarAx}, envisager les \'enonc\'es suivants :
  
 \begin{sugns}
 \label{sugKollar} 
 Soit $A$ un anneau de valuation discr\`ete de corps r\'esiduel   $F$.
  Soit $\X$ un $A$-sch\'ema r\'egulier, propre et plat sur $A$,
 \`a fibre g\'en\'erique $X$  lisse g\'eo\-m\'e\-tri\-quement connexe, \`a fibre sp\'eciale un diviseur $Y/F$ \`a croisements normaux
 stricts. Si $X/K$ est une vari\'et\'e  s\'eparablement rationnellement connexe, alors 
 
 (a) il  existe une composante r\'eduite $Y_{i}$ de $Y$ qui est g\'eo\-m\'e\-tri\-quement  int\`egre sur $F$;

(b) mieux, il existe une $F$-vari\'et\'e rationnellement connexe $Z$
et un $F$-morphisme de $Z$ dans $Y$;

(c) encore mieux, il existe une composante r\'eduite $Y_{i}$ de $Y$ qui est une 
 $F$-vari\'et\'e rationnellement connexe.
  \end{sugns}
  
  Dans cette direction, on a les r\'esultats suivants.

  \begin{thm}(Koll\'ar)\cite{KollarAx}
  \label{KollarAx}
   Soit $F$ un corps de caract\'eristique z\'ero.
  Soit $C$ une courbe lisse sur  $F$, soit
 $A$ l'anneau local de $C$ en un point ferm\'e de corps r\'esiduel $E$,
 soit $\X$ un $A$-sch\'ema r\'egulier, propre et plat sur $A$,
 \`a fibre g\'en\'erique $X$  lisse, \`a fibre sp\'eciale un diviseur $Y/E$ \`a croisements normaux
 stricts. Si $X/K$ est une vari\'et\'e de Fano, alors  il  existe une composante r\'eduite $Y_{i}$ de $Y$ qui est g\'eo\-m\'e\-tri\-quement irr\'eductible sur $E$.
  \end{thm}
  
 Toute hypersurface  est une d\'eg\'en\'erescence
 d'une hypersurface lisse de m\^eme degr\'e.
Le r\'esultat  de Koll\'ar \'etablit ainsi  la conjecture  \ref{Ax} (Ax)  en caract\'eristique z\'ero :
  toute $F$-hypersurface de degr\'e $d$ dans $\bP^n_{F}$ avec $n \geq d$ contient une
  sous-$F$-vari\'et\'e g\'eom\'etriquement int\`egre.

  \begin{cor}
  \label{grandcorpsKollar}
  Soit $k$ un corps de caract\'eristique z\'ero. Il existe un corps $L$
  contenant $k$ poss\'edant les propri\'et\'es suivantes :
  
  (i) Le corps $k$ est alg\'ebriquement ferm\'e dans $L$.
  
  (ii) Le corps $L$ est union de corps de fonctions de $k$-vari\'et\'es g\'eom\'etriquement int\`egres.
  
  (iii) Toute $L$-vari\'et\'e g\'eom\'etriquement int\`egre  poss\`ede un point $L$-rationnel
  (le corps $L$ est {\og pseudo-alg\'ebriquement clos \fg}).
  
  (iv) Le corps $L$ est un corps $C_{1}$.
   \end{cor}
   
   \begin{proof}  La construction du paragraphe \ref{rappels}  donne un corps $L$ satisfaisant
   les propri\'et\'es (i) \`a (iii).  Le point (iv) est alors une application de la conjecture d'Ax. 
      \end{proof}  
  
  \begin{thm}(Starr)\cite{Starr} 
   \label{starrgeomirr}
   Soit $F$ un corps parfait contenant un corps alg\'e\-bri\-quement clos.
  Soit $C$ une courbe lisse sur   $F$, soit
 $A$ l'anneau local de $C$ en un point ferm\'e de corps r\'esiduel $E$,
 soit $\X$ un $A$-sch\'ema r\'egulier, propre et plat sur $A$,
 \`a fibre g\'en\'erique $X/K$  lisse, \`a fibre sp\'eciale un diviseur $Y/E$ \`a croisements normaux
 stricts. Si $X $ est une $K$-vari\'et\'e s\'eparablement rationnellement connexe, alors  il  existe une composante r\'eduite $Y_{i}$ de $Y$ qui est g\'eo\-m\'e\-tri\-quement irr\'eductible sur $E$.
  \end{thm}

  Le th\'eor\`eme suivant  implique en particulier que le th\'eor\`eme \ref{KollarAx}
  vaut plus g\'en\'eralement lorsque la fibre g\'en\'erique est une vari\'et\'e rationnellement connexe.

\begin{thm}(Hogadi et Xu)\cite{HogadiXu}
\label{HogadiXu}
 Soient $F$ un corps de caract\'eristique z\'ero, 
   $C$ une $F$-courbe lisse,  
 $A$ l'anneau local de $C$ en un point ferm\'e $P$, et   $E$ le  corps r\'esiduel en $P$.
Soit $\X$ un $A$-sch\'ema 
  propre et plat sur $A$,
 de fibre g\'en\'erique $X$  une $K$-vari\'et\'e rationnellement connexe. Alors
 
 (a)  Il existe une $E$-vari\'et\'e rationnellement connexe $Z$
et un $E$-morphisme de $Z$ dans la fibre $Y/E$ de $\X \to C$ en $P$.
 
(b) Si  $\X$ est r\'egulier, connexe, de dimension relative au plus 3, et  si la fibre sp\'eciale est un diviseur $Y/E$ 
\`a croisements normaux  stricts,
alors il existe une composante r\'eduite $Y_{i}$ de $Y$  qui est une 
 $E$-vari\'et\'e rationnellement connexe.
\end{thm}

Sous l'hypoth\`ese suppl\'ementaire que $F$ contient un corps alg\'ebriquement clos,
le r\'esultat (a) avait \'et\'e \'etabli ant\'erieurement par de Jong.
 
  \begin{cor}
  \label{grandcorpsHogadiXu}
  Soit $k$ un corps de caract\'eristique z\'ero. Il existe un corps $L$
  contenant $k$ poss\'edant les propri\'et\'es suivantes :
  
  (i) Le corps $k$ est alg\'ebriquement ferm\'e dans $L$.
  
  (ii) Le corps $L$ est union de corps de fonctions de $k$-vari\'et\'es rationnellement connexes.
  
  (iii) Toute $L$-vari\'et\'e rationnellement connexe  poss\`ede un point $L$-rationnel.
   
   (iv) Le corps $L$ est un corps $C_{1}$.
   \end{cor}
   \begin{proof}  
    On reprend la construction du
paragraphe  \ref{rappels} mais \`a la place des $F$-vari\'et\'es g\'eom\'etriquement int\`egres
 on utilise les $F$-vari\'et\'es int\`egres $F$-birationnelles \`a une $F$-vari\'et\'e
 rationnellement connexe. La conditions (Stab) du paragraphe \ref{rappels} est satisfaite 
 gr\^ace au th\'eor\`eme   \ref{ratconsurratcon}
 (cons\'equence du th\'eor\`eme de Graber, Harris et Starr).
 Le corps $L$ ainsi construit satisfait les propri\'et\'es (i) \`a (iii).
  Toute hypersurface  est une d\'eg\'en\'erescence
 d'une hypersurface lisse de m\^eme degr\'e.
Le th\'eor\`eme    \ref{HogadiXu}  implique donc 
 que le corps $L$ est un corps $C_{1}$ 
 (voir \cite{HogadiXu}, Cor. 1.5). 
      \end{proof}

Une variante de la d\'emonstration du th\'eor\`eme \ref{SRCsurF((t))} 
permet de g\'en\'eraliser  la partie (a) du  th\'eor\`eme \ref{HogadiXu}.
\begin{thm}
\label{HogadiXugeneralise}
Soit $A$ un anneau de valuation discr\`ete, de corps des fractions $K$ et de corps r\'esiduel $F$ de caract\'eristique nulle. Soit $\X $ un $A$-sch\'ema  projectif  et plat sur $A$,
de  fibre g\'en\'erique une $K$-vari\'et\'e rationnellement connexe.
  Alors il existe une $F$-vari\'et\'e rationnellement connexe $Z$
et un $F$-morphisme de $Z$ dans la fibre sp\'eciale $Y=\X\times_{A}F$.
\end{thm}

 \begin{proof} 
Pour \'etablir le r\'esultat, on peut remplacer $A$ par son com\-pl\'et\'e.
Comme la caract\'eristique de $F$ est  nulle, ce compl\'et\'e est isomorphe \`a $F[[t]]$.
On est donc r\'eduit au cas $A=F[[t]]$.
Soit $R$ le hens\'elis\'e de $F[t]$ en $t=0$. 
 Soit $L$ son corps des fractions. On a $\hat{R}=A$ et $\hat{L}=K$.
La d\'emonstration du th\'eor\`eme
\ref{SRCsurF((t))}
montre qu'il existe   un $R$-sch\'ema   projectif et plat $\X_{1}$ (non n\'ecessairement r\'egulier) 
de fibre g\'en\'erique une $L$-vari\'et\'e  rationnellement connexe,
 tel que $\X_{1} \times_{R}F \simeq \X\times_{F[[t]]}F=Y$, en d'autres termes,
    la fibre sp\'eciale $Y$ de $\X_{1}$ est $F$-isomorphe \`a la fibre sp\'eciale de $\X$.
    D'apr\`es le th\'eor\`eme  \ref{HogadiXu}  (Hogadi et Xu), il existe
    une $F$-vari\'et\'e rationnellement connexe $Z$ et un $F$-morphisme $Z \to  \X_{1}\times_{R}F$. 
 \end{proof}

\begin{rmq} (Wittenberg)
Soit $A$ un anneau de valuation discr\`ete de corps des fractions $K$ et de corps
r\'esiduel $F$ de caract\'eristique z\'ero. Soit $\X$ un $A$-sch\'ema 
r\'egulier, propre et plat sur $A$,
 \`a fibre g\'en\'erique $X$  lisse, \`a fibre sp\'eciale un diviseur $Y/F$ \`a croisements normaux
 stricts. Lorsque la fibre g\'en\'erique $X$ est une $K$-compactification lisse d'un  espace homog\`ene d'un $K$-groupe alg\'ebrique lin\'eaire connexe, on peut facilement \'etablir le th\'eor\`eme \ref{KollarAx}
 et l'\'enonc\'e (a) du th\'eor\`eme \ref{HogadiXu}.
 On remplace $A$ par $F[[t]]$. Comme rappel\'e au paragraphe \ref{rappels},
   le corps $F$ est alg\'ebriquement ferm\'e dans un corps $E$ de dimension cohomologique
 $cd(E) \leq 1$, corps qui est union  de corps de fonctions de $F$-vari\'et\'es d'un type
 sp\'ecial, en particulier rationnellement connexes.
  Le corps $L$ limite inductive des corps $E((t^{1/n}))$ a le m\^eme groupe
 de Galois que $E$. Il est donc de dimension cohomologique 1, et 
 $X$ a un $L$-point. Ceci implique l'existence d'une $F$-application rationnelle 
 d'une $F$-vari\'et\'e rationnellement connexe  $Z$
 dans une composante r\'eduite de la fibre sp\'eciale, composante qui \'etant lisse
 doit en particulier \^etre g\'eo\-m\'e\-tri\-quement int\`egre.  
  \end{rmq}

\begin{rmq}
De m\^eme que l'on ne peut esp\'erer \'etendre le th\'eor\`eme \ref{fibSRCm1}  \`a d'autres classes  
de vari\'et\'es que celle des vari\'et\'es rationnellement connexes, de m\^eme il semble 
d\'eraisonnable d'esp\'erer une r\'eponse positive \`a la suggestion \ref{sugKollar} (a) pour d'autres
classes que celle des vari\'et\'es rationnellement connexes, par exemple pour les
vari\'et\'es projectives et lisses $X$ telles que  $H^{i}(X,O_{X})=0$ pour $i \geq 1$,
ou telles que le groupe de Chow  de $X$ r\'eduit des z\'ero-cycles sur tout corps alg\'ebriquement clos
soit nul (voir le paragraphe suivant).
Starr (communication priv\'ee) a donn\'e un exemple de surface d'Enriques sur
un corps $K(t)$
 telle que pour tout mod\`ele  de type (R) de cette surface sur l'anneau local de $K[t]$ en $t=0$,
  \`a croisements normaux stricts,
aucune composante r\'eduite ne soit g\'eom\'etriquement int\`egre.
 \end{rmq}
\section{Autour du th\'eor\`eme de Chevalley-Warning  : Vari\'et\'es dont le groupe de Chow g\'eom\'e\-trique 
   est trivial}

Le th\'eor\`eme de Chevalley-Warning a fait l'objet de plusieurs g\'en\'eralisations (Ax, Katz, Esnault, voir \cite{CL}).

\begin{thm} (Weil, 1954)   Toute surface projective lisse g\'eo\-m\'e\-tri\-quement rationnelle sur un corps fini poss\`ede un point rationnel. 
\end{thm}

\begin{thm} (formule de Woods Hole 1964, de Lefschetz-Verdier, voir Grothendieck/Illusie  SGA 5 III, Katz  SGA7 XXII)  
\label{woodshole}
Soit $\F$ un corps fini de caract\'eristique $p$.  Soit $X/\F$  une vari\'et\'e propre. Si $H^0(X,O_{X})=\F $ et si   $H^r(X,O_{X})=0$ pour $r\geq 1$, alors
le nombre de points rationnels de $X$ est congru \`a 1 modulo $p$.
 \end{thm}

En caract\'eristique nulle, les groupes   $H^r(X,O_{X}) \hskip1mm (r\geq 1)$ s'annulent pour une vari\'et\'e de Fano,
  mais on ne sait pas le d\'emontrer en caract\'eristique positive (sauf en dimension au plus 3, le cas de la dimension 3
 \'etant d\^u \`a Shepherd-Barron).

\medskip

H. Esnault a  obtenu le r\'esultat suivant.

\begin{thm} (Esnault 2003) \cite{Esnault1}
\label{Esnault1}  
Soit $\F$ un corps fini de cardinal $q$. Pour $X/\F$ lisse, projective, g\'eo\-m\'e\-tri\-quement int\`egre, 
et \ $\Omega$ un corps alg\'ebri\-quement clos contenant  le corps $\F (X)$, si l'on a $A_{0}(X_{\Omega})=0$, alors le nombre de points
$\F$-rationnels de $X$ est congru \`a 1 modulo $q$. 
\end{thm}

Soit $l$ un nombre premier, $l \neq {\rm car}(\F)$.
Par un argument remontant \`a Spencer Bloch  et d\'evelopp\'e par Bloch et Srinivas, 
l'hypoth\`ese assure que la cohomologie $l$-adique de $X$ est de coniveau 1, c'est-\`a-dire qu'elle satisfait
$H^{i}_{\et}({\overline X},\Q_{l})=N^1H^{i}_{\et}({\overline X},\Q_{l})$ pour tout $i \geq 1$ (toute classe de cohomologie s'annule sur un ouvert de Zariski non vide).
Sous cette condition, H. Esnault utilise des r\'esultats de Deligne pour \'etablir la congruence
annonc\'ee.

Ce th\'eor\`eme s'applique pour les vari\'et\'es rationnellement connexes par cha\^{\i}nes,
et en particulier pour les  vari\'et\'es de Fano, \`a la diff\'erence du th\'eor\`eme \ref{woodshole}.

A noter que le th\'eor\`eme s'applique aussi pour des vari\'et\'es qui ne sont pas rationnellement connexes,
comme les  surfaces d'Enriques et aussi  certaines surfaces de type g\'en\'eral.

\medskip

Comme dans l'\'enonc\'e initial de Chevalley-Warning et dans l'\'enonc\'e du th\'eor\`eme de Tsen,
on dispose de versions portant sur les fibres sp\'eciales, singuli\`eres, de telles vari\'et\'es.

\begin{thm} (N. Fakhruddin et C.S. Rajan, 2004)   \cite{FR}  
Soit $f : X \to Y$ un morphisme propre dominant de vari\'et\'es lisses
 et g\'eo\-m\'e\-tri\-quement irr\'eductibles sur un corps fini $\F$ de cardinal $q$. Soit  $Z$ la fibre g\'en\'erique,
 suppos\'ee g\'eo\-m\'e\-tri\-quement int\`egre. 
 Soit \ $ \overline{\F(Y)} $ une cl\^oture alg\'ebrique du corps des fonctions $\F(Y)$.
 Si l'on a 
 $A_{0}(Z_ {\overline{\F(Y)}})  =0$,
 alors  pour tout point   $y \in Y(\F )$, le cardinal de  $X_{y}(\F)$ est congru \`a 1 modulo $q$.
 Si l'hypoth\`ese $X$ lisse est omise mais si  la fibre g\'en\'erique $Z$ est lisse,
on a  $X_{y}(\F)\neq \emptyset$ pour tout $y \in Y(\F)$.
\end{thm}

Donc sur toute d\'eg\'en\'erescence de vari\'et\'e RCC (lisse) il y a un $\F$-point.
Ceci vaut aussi sur une d\'eg\'en\'erescence d'une surface d'Enriques
ou de    certaines surfaces de type g\'en\'eral.
\begin{thm} (Esnault  \cite{Esnault2} \cite{Esnault3}; Esnault et Xu \cite{Esnault4})
Soit $A$ un anneau de valuation discr\`ete complet de corps des fractions $K$ et de corps r\'esiduel
$\F$ fini de cardinal $q$. Soit ${\X}$ un $A$-sch\'ema int\`egre   propre et plat. 
Soit  $l$ un nombre premier, $l \neq {\rm car}(\F)$.
Supposons la fibre g\'en\'erique g\'eo\-m\'e\-tri\-quement int\`egre, lisse et  
 \`a cohomologie $l$-adique de coniveau 1. Soit $Y/\F$ la fibre sp\'eciale.  Alors
  
(i)   $Y(\F) \neq \emptyset$;

(ii)  si de plus  $\X$ est r\'egulier, alors ${\rm card} (Y(\F)) \equiv 1 \hskip1mm {\rm mod} \hskip1mm q$.
\end{thm}

  L'hypoth\`ese sur la cohomologie est satisfaite si   $A_{0} (X\times_{K}{\Omega})=0$, o\`u $\Omega$
  est un corps alg\'ebriquement clos contenant $K(X)$, en particulier pour les vari\'et\'es RCC mais aussi pour les surfaces d'Enriques et certaines surfaces de type g\'en\'eral.

  \medskip
 
 En particulier il y a un point rationnel sur la fibre sp\'eciale. En particulier si toutes les composantes de la fibre sp\'eciale sont lisses, alors l'une d'entre elles est g\'eo\-m\'e\-tri\-quement int\`egre sur $\F_{q}$.
 
 \bigskip
 
 Il y a des th\'eor\`emes de g\'eom\'etrie alg\'ebrique qui se d\'emontrent par r\'eduction au cas
 des corps finis. On part d'une vari\'et\'e sur un corps $k$. Une telle vari\'et\'e
 est obtenue par changement de  base $A \to k$
 \`a partir d'un $A$-sch\'ema de type fini, pour une $\Z$-alg\`ebre de type fini $A$ convenable.
 On r\'eduit ensuite aux points ferm\'es de $A$ (leurs corps r\'esiduels sont finis)
 et on applique les r\'esultats obtenus sur les corps finis.

 Les th\'eor\`emes \'etablis par H. Esnault sont de ce point de vue {\og trop bons \fg} :
 la classe des $K$-vari\'et\'es auxquelles  ses r\'esultats s'appliquent est plus large que celle
 des $K$-vari\'et\'es rationnellement connexes.
  On ne peut donc esp\'erer les utiliser pour
 \'etablir des r\'esultats comme le th\'eor\`eme \ref{GHSJ} (Graber-Harris-Starr) ou le th\'eor\`eme   \ref{KollarAx} ci-dessus (Koll\'ar) --  pas plus
 d'ailleurs que l'on ne pouvait utiliser le th\'eor\`eme de Chevalley-Warning pour
 \'etablir le th\'eor\`eme de Tsen ou la conjecture d'Ax. Un obstacle essentiel semble \^etre le fait bien
 connu suivant : il existe des polyn\^omes   en une variable sur $\Z$
 qui n'ont pas de z\'ero sur $\Q$ mais dont la r\'eduction en tout premier $p$ sauf un nombre fini
 a un z\'ero, par exemple $(x^2-a)(x^2-b)(x^2-ab)$, avec $a, b \in \Z$ non carr\'es.
 
 Les r\'esultats sur les corps finis peuvent n\'eanmoins en sugg\'erer d'autres
 sur les corps de fonctions d'une variable. On en trouvera un exemple r\'ecent
 dans \cite{Kahn}, \S 9.8, Remarque 3.
   
 \section{Approximation faible pour les vari\'et\'es rationnellement connexes}

\begin{sugn}
\label{appfaibleratcon}
  Soit $K$ un corps de fonctions d'une variable sur un corps alg\'ebriquement clos.
      Pour toute vari\'et\'e rationnellement connexe $X$ sur $K$, l'approximation faible vaut :
      pour tout ensemble fini $I$ de places $v$ de $K$, l'application diagonale
      $$X(K) \to \prod_{v \in I} X(K_{v})$$
      a une image dense. Ici $K_{v}$ est le compl\'et\'e de $K$ en $v$
      et $X(K_{v})$ est muni de la topologie induite par la topologie de la valuation sur  $K_{v}$.
\end{sugn}

     Des arguments \'el\'ementaires  (\cite{CTGille})  permettent d'\'etablir l'approximation faible en tout ensemble fini de places pour les compactifications lisses d'espaces homog\`enes de groupes lin\'eaires connexes, puis  pour
    les vari\'et\'es obtenues par fibrations en  de telles vari\'et\'es. On traite ainsi les intersections compl\`etes  lisses de deux quadriques dans $\bP^n$ pour $n \geq 4$.

\bigskip

\begin{thm} (Hassett-Tschinkel)\cite{HT1}
 Soit $K$ un corps de fonctions d'une variable sur un corps alg\'ebriquement clos de caract\'eristique z\'ero. Soit $X/K$ une $K$-vari\'et\'e  
   rationnellement connexe.
    Si $I$ est un ensemble fini de places  de $K$ de bonne r\'eduction pour $X/K$, alors 
    l'approximation faible vaut pour $X$ en ces places : l'application diagonale
     $X(K) \to \prod_{v \in I} X(K_{v})$ a une image dense.
     \end{thm}
    
    Ceci g\'en\'eralise un r\'esultat de Koll\'ar, Miyaoka et  Mori (cas on l'on demande une r\'eduction fix\'ee, sans obtenir d'approximations aux jets d'ordre sup\'erieur).
    
Le cas particulier des surfaces cubiques lisses avait \'et\'e trait\'e par Madore \cite{MadoreSMF}.

    Hassett et Tschinkel \cite{HT2} ont aussi des r\'esultats d'approximation en des places de mauvaise,
    mais pas trop mauvaise r\'eduction. Mais comme ces auteurs le notent,  le cas suivant est ouvert.

    \begin{qn}
    L'approximation faible en la  place $\lambda=0$ vaut-elle pour 
      la surface cubique
     $x^3+y^3+z^3+\lambda t^3=0$ sur le corps $K=\C(\lambda)$ ?
    \end{qn}

  Lorsque le nombre de variables est  suffisamment grand par rapport au degr\'e, on a pu \'etablir
 l'approximation faible en toutes les places. Voir la section \ref{apfaiblepartout}  ci-dessous.

\section{$R$-\'equivalence sur les vari\'et\'es  rationnellement connexes}
\label{Requivratcon}

Soient $k$ un corps non alg\'ebriquement clos et $X$ une $k$-vari\'et\'e (s\'eparablement) rationnellement connexe.
Que sait-on sur l'ensemble $X(k)/R$ ?

\medskip

\begin{thm}(Koll\'ar)\cite{Kollarannmath}
 Soit $K$ un corps local usuel (localement compact) et soit $X$ une $K$-vari\'et\'e 
 s\'eparablement rationnellement connexe.
Alors la $R$-\'equivalence sur $X(K)$ est une relation ouverte. L'ensemble $X(K)/R$ est fini.
Dans le cas $K=\R$ les classes de $R$-\'equivalence co\"{\i}ncident avec les
composantes connexes de $X(\R)$.
\end{thm}

Ce r\'esultat est une vaste g\'en\'eralisation de cas particuliers ant\'erieurement connus
(surfaces fibr\'ees en coniques, compactifications de groupes alg\'ebriques lin\'eaires connexes,
hypersurfaces cubiques lisses,
 intersections lisses de deux quadriques dans $\bP^n$ pour $nÊ\geq 4$).

\begin{sugn} (Koll\'ar)
Soient $\F$ un corps fini et $X$ une $\F$-vari\'et\'e 
 s\'epa\-ra\-blement rationnellement
connexe. Alors tous les points de $X(\F)$ sont $R$-\'equivalents : l'ensemble $X(\F)/R$ 
a un \'el\'ement.
 \end{sugn}

 Swinnerton-Dyer montra qu'il en est ainsi pour les surfaces cubiques lisses.
 Ce r\'esultat a \'et\'e  r\'ecemment \'etendu par J. Koll\'ar \cite{DK} \`a toutes les hypersurfaces cubiques lisses
 sur un corps fini de cardinal au moins 8.

\begin{thm} (Koll\'ar-Szab\'o)\cite{KollarSzabo}
  Soient $\F$ un corps fini et $X$ une $\F$-vari\'et\'e 
     s\'eparablement rationnellement
connexe. Si l'ordre de $\F$ est plus grand qu'une certaine constante qui d\'epend seulement de
la g\'eom\'etrie de $X$ alors $X(\F)/R$ est r\'eduit \`a un point.
\end{thm}

\begin{thm} (Koll\'ar-Szab\'o)\cite{KollarSzabo} Soit $K$ un corps local non archim\'edien de corps r\'esiduel le corps fini $\F$. Soit $A$
l'anneau de la valuation.  Soit $\X$ un $A$-sch\'ema r\'egulier, int\`egre, projectif  et
plat sur $A$, de fibre sp\'eciale $Y/\F$ une $\F$-vari\'et\'e 
 s\'eparablement rationnellement
connexe -- ce qui implique que la fibre g\'en\'erique $X=\X\times_{A}K$ est SRC.
Si l'ordre de $\F$ est plus grand qu'une certaine constante qui d\'epend seulement de
la g\'eom\'etrie de $X$ alors $X(K)/R$ est r\'eduit \`a un point.
\end{thm}

Ici encore on se demande si la condition sur l'ordre du corps r\'esiduel est n\'ecessaire.
A tout le moins, le r\'esultat ci-dessus implique  :

\begin{thm} \cite{KollarSzabo}
 Soient $K$ un corps de nombres et $X/K$ une $K$-vari\'et\'e 
  rationnellement
connexe. Alors pour presque toute place $v$ de $K$, notant $K_{v}$ le compl\'et\'e de
$K$ en $v$, on a  ${\rm card} \hskip1mm  X(K_{v})/R=1$.
\end{thm}

\medskip

Soit $A$ un anneau de valuation discr\`ete de corps des fractions $K$ et de corps r\'esiduel $F$.
Soit $\X$ un $A$-sch\'ema   int\`egre, propre et lisse. Soit $X=\X\times_{A}K$
la fibre g\'en\'erique et $Y=\X\times_{A}F$ la fibre sp\'eciale.

La sp\'ecialisation $X(K)=\X(A) \to Y(F)$ passe au quotient par la
$R$-\'equivalence (voir  \cite{Madorethese}). On a donc
  une application  de sp\'ecialisation :
$$X(K)/R \to Y(F)/R.$$

\begin{thm} (Koll\'ar) \cite{Kollarjap} Dans la situation ci-dessus, si $Y/F$ est SRC, et si $A$ est 
 hens\'elien, alors
l'application de sp\'ecialisation $X(K)/R \to Y(F)/R$ est une bijection.
\end{thm}

\medskip

On dit qu'un corps $K$ est fertile (les anglo-saxons disent {\og large field  \fg})  
si sur toute $K$-vari\'et\'e lisse int\`egre avec un $K$-point
les $K$-points sont denses pour la topologie de Zariski.  

Exemples : 

(a) Une extension alg\'ebrique infinie d'un corps fini (estimations de Lang--Weil).

 (b) Un corps    local usuel (non archim\'edien, \`a corps r\'esiduel fini), plus g\'en\'eralement le corps des fractions d'un anneau de valuation discr\`ete hens\'elien de corps r\'esiduel quelconque.
 
(c) Le corps  $\R$ des r\'eels, plus g\'en\'eralement un corps r\'eel clos, plus g\'en\'e\-ra\-lement
un corps dont le groupe de Galois absolu est un pro-$p$-groupe ($p$ \'etant un nombre premier).

(d) Un corps pseudo-alg\'ebriquement clos.

\begin{thm} (Koll\'ar)  Soient $K$ un corps fertile  et $X$ une $K$-vari\'et\'e 
 s\'eparablement  rationnellement connexe.

(1) \cite{Kollarannmath}  Pour tout point $M \in X(K)$, il existe un $K$-morphisme tr\`es libre $f : \bP^1_{K} \to X$
tel que $M$ appartienne \`a $f(\bP^1(K))$.

(2) \cite{Kollarjap}  Si deux $K$-points sont $R$-\'equivalents, alors il existe un $K$-mor\-phisme $\bP^1_{K} \to X$
tel que ces deux points soient dans l'image de $\bP^1(K)$.
\end{thm}

\begin{cor}  \cite{Kollarjap}
 Pour $K$ corps fertile et $X/K$ comme ci-dessus, pour tout ouvert de Zariski non vide $U \subset X$,
l'application $U(K)/R \to X(K)/R$ est bijective.
\end{cor}

On ne sait pas si les deux \'enonc\'es pr\'ec\'edents valent  sur un corps
$K$ infini quelconque.

\begin{thm} (Koll\'ar)  \cite{Kollarjap}  Soit $K$ un corps local usuel, soit $f : X \to Y$ un $K$-morphisme projectif et lisse
de $K$-vari\'et\'es lisses, dont les fibres g\'eom\'etriques 
  sont des vari\'et\'es SRC. L'application $Y(K) \to \N$
qui \`a un point $y \in Y(K)$ associe le cardinal de $X_{y}(K)/R$ est semi-continue sup\'erieurement :
tout point de $Y(K)$ admet un voisinage (pour la topologie sur $Y(K)$ d\'efinie par celle du corps local $K$)  tel que pour $z$ dans ce voisinage le cardinal de
$X_{z}(K)/R$ soit au plus \'egal \`a celui de $X_{y}(K)/R$.
\end{thm}

\begin{qn}
  \cite{Kollarjap}  Le cardinal de $X_{y}(K)/R$ est-il localement constant  quand $y$ varie dans $Y(K)$ ?
  \end{qn}

\begin{qn} 
\label{Requivcd1}
Soient $k$ un corps  et $X$ une $k$-vari\'et\'e
s\'eparablement rationnellement connexe.
Dans chacun des cas suivants :

(a)  $k=\C(C)$ est un corps de fonctions d'une variable sur les complexes,  

(b) $k=\C((t))$ est un corps de s\'eries formelles en une variable,

(c) $k$ est un corps  $C_{1}$,

(d) $k$ est un corps  parfait de dimension cohomologique $cd(k) \leq 1$,
 
 \noindent l'ensemble $X(k)/R$ a-t-il au plus un \'el\'ement  ?
 \label{Requivtriv}
\end{qn}

 On ne s'attend pas \`a une r\'eponse positive. 
  Cependant, 
pour $k$ de carac\-t\'eristique nulle, 
 sous la simple hypoth\`ese $cd(k) \leq 1$,
  c'est connu dans les cas suivants :

(i)  $X$ est une compactification lisse d'un groupe lin\'eaire connexe (\cite{RET}).

(ii)  $X$ est une surface fibr\'ee en coniques de degr\'e 4 sur la droite projective
(\cite{CTSkMan}).

(iii) $X$ est une intersection lisse de deux quadriques dans $\bP^n_{k}$
et $n \geq  5$  (\cite{CTSaSwD}, Thm. 3.27 (ii)).
 
(iv) Le corps $k$ est $C_{1}$,  la vari\'et\'e $X$ est une hypersurface cubique lisse dans $\bP^n_{k}$  
avec $n \geq 5$
(\cite{MadoreJNT}).

On a aussi le  r\'esultat  suivant, portant sur des vari\'et\'es singuli\`eres :

(v) Soit $k$ un corps de caract\'eristique nulle
  tel que toute forme quadratique sur $k$ en 3 variables ait un z\'ero non trivial.
Alors pour {\it toute}   surface cubique {\it singuli\`ere} $X \subset \bP^3_{k}$ 
l'ensemble $X(k)/R$ a au plus un \'el\'ement.   

Lorsque $X$ poss\`ede un point singulier $k$-rationnel, ceci est \'etabli dans \cite{MadoreJNT}, \S 1.
Dans le cas g\'en\'eral, on \'etablit ce r\'esultat en utilisant la classification des surfaces
cubiques singuli\`eres. Le seul cas non couvert par
les arguments donn\'es au \S 5 de  \cite{Madoremanmat} (voir aussi \cite{MadoreJNT}, Remarque 1)
est le cas des surfaces de Ch\^atelet (cas 7 p. 182 de  \cite{Madoremanmat}). Le r\'esultat dans ce cas
s'obtient en combinant le Th\'eor\`eme 8.6 (d) de  \cite{CTSaSwD}  et les r\'esultats de \cite{coraytsfasman}. 

\medskip

Une r\'eponse positive \`a la question \ref{Requivtriv}  pour les surfaces (projectives et lisses) g\'eom\'etriquement ration\-nelles d\'efinies sur $\C(t)$ impliquerait l'unirationalit\'e des vari\'et\'es de dimen\-sion 3 sur $\C$ qui admettent une fibration en coniques sur le plan projectif. Il s'agit l\`a d'une question largement ouverte.

\begin{qn} 
Soient $K$ un corps de nombres et  $X$  une $K$-vari\'et\'e 
 rationnellement connexe.
Le quotient $X(K)/R$ est-il fini ?
\end{qn}

 C'est connu  dans les cas suivants :

(i)  La vari\'et\'e $X$ est une compactification lisse d'un groupe lin\'eaire connexe $G$. 
L'immersion ouverte $G \subset X$ induit une bijection $G(k)/R \simeq X(k)/R$
(\cite{GilleTAMS}).  La finitude dans le cas g\'en\'eral est due \`a Gille \cite{GilleIHES},
elle s'appuie sur des r\'esultats ant\'erieurs de Margulis (groupes semi-simples simplement connexes)
et CT-Sansuc (\cite{RET}, cas des tores alg\'ebriques).

(ii) La vari\'et\'e  $X$ est une surface fibr\'ee  en coniques de degr\'e 4 sur la droite projective
(CT-Sansuc, cf. \cite{CTSkMan}).

(iii) La vari\'et\'e $X$ est une intersection lisse de deux quadriques dans $\bP^n_{K}$
et $n \geq 6$ (\cite{CTSaSwD}).

 \medskip

La question de la finitude de $X(K)/R$ sur $K$ un corps de nombres est ouverte pour les compactifications lisses d'espaces homog\`enes de groupes lin\'eaires connexes, m\^eme en supposant les groupes d'isotropie g\'eom\'etrique connexes.

\medskip

On pourrait se poser la question de la finitude de $X(K)/R$ pour $X/K$ rationnellement connexe
et $K$ de type fini sur l'un quelconque des corps suivants :  un corps fini,  $\Q$,   $\C$,   $\R$,   $\Q_{p}$.

On a par exemple la finitude dans ce cadre dans le cas (ii) ci-dessus (\cite{CTSkMan}), et c'est 
une question ouverte lorsque $X$
est de dimension 2, i.e. est une surface g\'eo\-m\'e\-tri\-quement rationnelle.
Dans le cas (i), on a la finitude lorsque $G$ est 
 un tore \cite{RET}.   
 C'est une question largement ouverte pour $G$ un groupe  lin\'eaire quelconque.
  
Mais, sur chacun des   corps   $\Q(t)$, $\R(t)$, $\R((t))$,
la r\'eunion pour tout $n \geq 1$ des $\R((t^{1/n}))$ (qui est un corps r\'eel clos non archim\'edien),
Koll\'ar \cite{Kollarjap}  a construit des exemples d'hypersurfaces lisses $X$ de degr\'e 4 dans $\bP^n_{K}$, avec $n$ arbitrairement grand, telles que $X(K)/R$ soit infini.

\section{\'Equivalence rationnelle sur les z\'ero-cycles des vari\'et\'es rationnellement connexes}

\begin{prop} 
\label{CTinv06}
\cite{CTinvmath06} Soient $k$ un corps et $X$ une $k$-vari\'et\'e 
 RCC. Il existe un entier $N=N(X)>0$ tel que pour toute extension de corps $L/k$ on ait
  $NA_{0}(X\times_{k}L)=0$.
\end{prop}

Soient $\F$ un corps fini et $X$ une $\F$-vari\'et\'e s\'eparablement rationnellement connexe.
Du th\'eor\`eme de Koll\'ar et Szab\'o  \cite{KollarSzabo} il r\'esulte que l'on a $A_{0}(X)=0$.
Mais ceci n'est qu'un cas particulier d'un th\'eor\`eme g\'en\'eral en th\'eorie du corps de classes sup\'erieur :

\begin{thm}  (K. Kato et S. Saito, 1983)
 Soient $\F$ un corps fini et $X$ une $\F$-vari\'et\'e projective et lisse g\'eo\-m\'e\-tri\-quement int\`egre.
Soit $Alb_{X}$ la vari\'et\'e d'Albanese de $X$ (c'est une vari\'et\'e ab\'elienne) et $\mu$ le $\F$-groupe fini commutatif dual de la torsion du
groupe de N\'eron-Severi g\'eom\'etrique de $X$. Le groupe $A_{0}(X)$ est fini, et l'on a une suite exacte
$$ 0 \to H^1(\F,\mu) \to A_{0}(X)  \to Alb_{X}(\F) \to 0.$$
\end{thm}

\begin{qn}
Soient $k$ un corps  et $X$ une $k$-vari\'et\'e
s\'eparablement rationnellement connexe.
Dans chacun des cas suivants :

(a)  $k=\C(C)$ est un corps de fonctions d'une variable sur les complexes, 

(b) $k=\C((t))$ est un corps de s\'eries formelles en une variable,

(c) $k$ est un corps parfait de dimension coho\-mo\-lo\-gique~1,
 
\noindent a-t-on  $A_{0}(X)=0$ ?

\end{qn}

On ne s'attend pas \`a une r\'eponse positive. Cependant, 
pour $k$ de caract\'eristique nulle, sous la simple hypoth\`ese $cd(k) \leq 1$, il en est ainsi dans chacun des cas suivants :

(i) Compactification lisse d'espace  homog\`ene principal de groupe alg\'ebrique lin\'eaire connexe
  (\cite{RET}).

(ii) Surface  SRC, i.e. surface  g\'eo\-m\'e\-tri\-quement rationnelle. La situation est ici bien meilleure 
que pour la $R$-\'equivalence (voir la question \ref{Requivcd1}). On \'etablit  $A_{0}(X)=0$ par des m\'ethodes de  K-th\'eorie alg\'ebrique (\cite{CTinvmath83}).

 (iii) Hypersurface cubique lisse dans $\bP^n_{k}$ ($n \geq 3$) avec un $k$-point, pour $n \geq 3$.
 Soit $P \in X(k)$. Pour \'etablir (iii), il suffit de montrer que   tout $k$-point $M$ est rationnellement \'equivalent au point $P$ (on applique ensuite cet \'enonc\'e sur toute extension finie de $k$.)
 Soit $L \subset  \bP^n_{k}$ un espace lin\'eaire de dimension 3 contenant  $P$ et $M$.
 Soit $Y=X \cap L \subset  L \simeq \bP^3_{k}$ la surface cubique d\'ecoup\'ee par $L$.
 Si $Y$ est singuli\`ere, alors $P$ et $M$ sont $R$-\'equivalents sur $Y$, donc sur $X$ :
 voir l'\'enonc\'e (v) apr\`es la question  \ref{Requivtriv}. Si $Y$ est non singuli\`ere, on
 a $A_{0}(Y)=0$  d'apr\`es le point (ii) ci-dessus. Dans tous les cas on voit que $P$ et $M$
 sont rationnellement \'equivalents.

(iv) Intersection lisse de deux quadriques dans $\bP^n_{k}$  avec un $k$-point,  
pour $nÊ\geq 5$. 
Ceci r\'esulte de l'\'enonc\'e (iii) suivant la question \ref{Requivtriv}. Une adaptation 
de l'argument donn\'e ci-dessus pour les hypersurfaces cubiques
devrait donner le r\'esultat  pour $nÊ\geq 4$.

\begin{thm} (CT-Ischebeck 1981)
  Soit $X$ une $\R$-vari\'et\'e projective et lisse g\'eo\-m\'e\-tri\-quement int\`egre
avec $X(\R) \neq \emptyset$.  Soit $s $ le nombre de composantes connexes
de $X(\R)$.

Le sous-groupe $2A_{0}(X)$ est le sous-groupe divisible maximal de $A_{0}(X)$
et le quotient $A_{0}(X)/2A_{0}(X)=(\Z/2)^{s-1}$.

En particulier si $X$ est rationnellement connexe et $X(\R)\neq \emptyset$, alors 
$A_{0}(X)$ est fini et
$A_{0}(X)=(\Z/2)^{s-1}$.
\end{thm}
 
Soit $R$ un corps r\'eel clos.  
Knebusch et Delfs ont montr\'e comment l'on peut, pour toute $R$-vari\'et\'e alg\'ebrique $X$,
donner une d\'efinition ad\'equate des {\og composantes connexes \fg} de $X(R)$. Celles-ci sont en nombre fini. Le th\'eor\`eme ci-dessus vaut dans ce cadre plus large. On comparera ceci avec la remarque finale de la section \ref{Requivratcon}.

 \begin{qn}
Soient $K$ un corps $p$-adique (extension finie de $\Q_{p}$) 
et $X$ une $K$-vari\'et\'e rationnellement connexe.
Le groupe $A_{0}(X)$ est-il fini ?
\end{qn}

Soit $A$ l'anneau de la valuation du corps local $K$, soit $\F$ le corps fini  r\'esiduel.  Voici des r\'esultats obtenus
dans cette direction.

(i) Si ${\rm dim}(X)=2$, le groupe $A_{0}(X)$ est fini (\cite{CTinvmath83}).

(ii) Si $X$ est une intersection lisse de deux quadriques dans $\bP^n_{K}, n \geq 4$ et $X(K) \neq \emptyset$,
le groupe $A_{0}(X)$ est fini   (\cite{CTSaSwD},  \cite{CTSko} et  \cite{ParimalaSuresh}).

(iii) Si $X$ est un fibr\'e en quadriques de dimension relative au moins 1 sur la droite projective, le groupe $A_{0}(X)$ est fini
(\cite{CTSko, ParimalaSuresh}).

(iv)   Si $X$ est une $K$-compactification lisse d'un $K$-groupe lin\'eaire connexe, alors $A_{0}(X)$ est somme d'un groupe fini et d'un groupe de torsion $p$-primaire (d'exposant fini) (\cite{CTinvmath06}).

(v)  (Koll\'ar-Szab\'o) \cite{KollarSzabo} Si $X$ a bonne r\'eduction SRC, i.e. s'il existe un $A$-sch\'ema $\X$ r\'egulier, int\`egre, propre et lisse
  de fibre sp\'eciale $Y/\F$ SRC, alors  $A_{0}(X)=0$.

 (vi) (S. Saito et K. Sato) \cite{SS07} Soit $X$ une $K$-vari\'et\'e projective, lisse, g\'eom\'etriquement connexe. Supposons que $X/K$ poss\`ede un mod\`ele
 $\X/A$  r\'egulier int\`egre, propre et plat,  de
  fibre sp\'eciale r\'eduite $Y_{{red}}/\F$ \`a croisements normaux stricts.
Alors le groupe $A_{0}(X)$ est somme directe d'un groupe fini et d'un groupe divisible par tout
entier premier \`a $p$.
 Si en outre  $X$ est une vari\'et\'e rationnellement connexe, alors
$A_{0}(X)$ est  somme d'un groupe fini et d'un groupe de torsion $p$-primaire d'exposant fini.

\medskip

On s'est longtemps pos\'e la question de savoir si pour toute vari\'et\'e projective
lisse $X$ sur un corps $p$-adique le sous-groupe de torsion de $A_{0}(X)$ est fini.
M. Asakura et S. Saito ont montr\'e r\'ecemment qu'il n'en est rien (exemples : surfaces
de degr\'e $d \geq5$ suffisamment g\'en\'erales dans $\bP^3$).

\medskip

 \begin{qn}
Soient $K$ un corps de type fini sur le corps premier
et $X$ une $K$-vari\'et\'e rationnellement connexe.
Le groupe $A_{0}(X)$ est-il fini ?
\end{qn}

C'est connu lorsque ${\rm dim}(X)=2$  et $X(K) \neq \emptyset$ (\cite{CTinvmath83}),
et lorsque $X$ est une compactification lisse d'un $K$-tore de dimension 3 (Merkur'ev, 2008).
 
 Mais  le cas g\'en\'eral des compactifications lisses de  tores sur un corps de nombres 
 est ouvert.
 
 De fa\c con g\'en\'erale, on se demande si pour toute vari\'et\'e $X$ connexe, projective et lisse   
 sur un corps $K$ de type fini sur le corps premier, le groupe $A_{0}(X)$ est un groupe de type fini.

\section{Vers les vari\'et\'es  sup\'erieurement rationnellement connexes}

\subsection{Deux exemples}

\subsubsection{Formes tordues d'hyperquadriques}
\label{extao}

  D. Tao \cite{Tao}  a obtenu les r\'esultats suivants.
  Soit $K$ un corps poss\'edant une alg\`ebre simple centrale $A$
de degr\'e $2n \geq 6$ dont la classe $[A]$ dans le groupe de Brauer de $K$ est non nulle et d'exposant 2. La condition $2.[A]=0$ assure l'existence sur $A$ d'une involution de premi\`ere esp\`ece $\sigma$
qu'on peut choisir orthogonale. A une telle situation est alors associ\'ee une $K$-vari\'et\'e $X$
qui est une forme tordue d'une quadrique lisse dans $\bP^{2n-1}$ et pour laquelle $$ {\rm Ker} ( \Br K \to \Br K(X) ) = \Z/2 = \Z.[A] \subset \Br K.$$
Il y a donc une obstruction \'el\'ementaire \`a l'existence d'un $K$-point,  
en particulier $X(K)=\emptyset$.

On peut trouver des alg\`ebres $A$ du type requis sur l'un quelconque des corps suivants :
 un corps $p$-adique, le corps des s\'eries formelles it\'er\'ees $\C((u))((v))$, un corps de
fonctions de deux variables sur $\C$.

Sur $K$ l'un quelconque de ces corps, pour $m \geq 4$, les quadriques dans $\bP^m_{K}$  ont un point rationnel.  Mais pour $m$ impair  les formes tordues obtenues n'ont pas de point $K$-rationnel.

 Si l'on consid\`ere une telle forme tordue $X$ sur le corps $K=\C((u))((v))$,  pour laquelle
 l'application  $\Br K \to \Br K(X)$ n'est pas injective, il r\'esulte de la proposition
 \ref{injbrfibres} que  la fibre sp\'eciale sur $F=\C((u))$ d'un mod\`ele propre, plat, r\'egulier de $X$  sur l'anneau de valuation discr\`ete $A=\C((u))[[v]]$ n'a aucune composante g\'eo\-m\'e\-tri\-quement int\`egre de multiplicit\'e  1.

\subsubsection{Une hypersurface cubique}
\label{exdj}

L'hypersurface cubique diagonale  $X \subset \bP^8_{K}$ de coefficients 
$$(1,u,u^2,v,vu,vu^2,v^2,v^2u,v^2u^2)$$
sur le corps  $K=\C((u))((v))$ n'a pas de point rationnel.
  
 La condition d'injectivit\'e sur le groupe
de Brauer $\Br K \hookrightarrow \Br K(X)$
est ici satisfaite, plus g\'en\'eralement, l'obstruction \'el\'ementaire s'annule : il en est
ainsi pour toute hypersurface lisse de dimension au moins 3 (cf. \cite{BoCTSk}).

L'anneau de valuation discr\`ete $A=\C((u))[[v]]$ a pour corps des fractions $K$.
La $K$-hypersurface cubique
$X$ admet un mod\`ele r\'egulier $\X$ projectif sur $A$ dont une composante
r\'eduite de la fibre sur $F=\C((u))$ est g\'eo\-m\'e\-tri\-quement int\`egre et de multiplicit\'e  1 :
 un ouvert est donn\'e par un ouvert du c\^one de $\bP^8_{F}$ d'\'equation homog\`ene $x^3+uy^3+u^2z^3=0$. 
  Mais cette composante
n'est pas rationnellement connexe, elle ne  poss\`ede m\^eme pas de $\C((u))$-point lisse.
De fait, la fibre sp\'eciale $Y$ de $\X$ ne saurait poss\'eder une composante g\'eom\'etriquement int\`egre rationnellement connexe de multiplicit\'e  1~ : d'apr\`es le th\'eor\`eme \ref{SRCsurF((t))}
toute telle composante 
poss\'ederait
des points lisses sur $\C((u))$, points qui seraient Zariski-denses car $\C((u))$ est fertile,
et l'on pourrait relever un $\C((u))$-point non situ\'e sur les autres composantes en un $K$-point
de $X$. Le m\^eme argument montre qu'aucune composante de multiplicit\'e 1 de la fibre sp\'eciale
d'un mod\`ele   $\X$  de type (R) de $X$,   \`a croisements normaux stricts,  n'est le but d'une application rationnelle
depuis une $\C((u))$-vari\'et\'e rationnellement connexe.

\subsection{Fibres sp\'eciales avec une composante g\'eom\'e\-tri\-que\-ment int\`egre de multiplicit\'e 1}

 Soit $A$ un anneau de valuation discr\`ete, $K$ son corps des fractions, $F$ son corps r\'esiduel.
Soit $\pi$ une uniformisante de $A$. Soit $\X/A$ un $A$-sch\'ema de type (R) (voir le paragraphe \ref{rappels}),
$X/K$ sa fibre g\'en\'erique, $Y/F$ sa fibre sp\'eciale. Si l'on a $X(K) \neq \emptyset$ alors
l'on a $\X(A)\neq \emptyset$. Comme $\X$ est r\'egulier, une $A$-section de $X/A$
rencontre $Y$ en un $F$-point $M$ poss\'edant les propri\'et\'es suivantes : il
est sur une unique composante r\'eduite de $Y$, il est lisse sur cette composante,
cette composante est de multiplicit\'e 1 et g\'eom\'etriquement int\`egre.
 Inversement, si $A$ est hens\'elien, un tel point $M$ se rel\`eve en un $K$-point de $X$.

On voit donc qu'une condition n\'ecessaire pour l'existence d'un $K$-point sur $X$
est l'existence d'une composante g\'eom\'etriquement int\`egre de multiplicit\'e 1
de la fibre sp\'eciale $Y$. Au paragraphe \ref{compgeom1} on a discut\'e cette propri\'et\'e.
 Par analogie avec les suggestions \ref{sugKollar} 
on est amen\'e ici \`a s'int\'eresser aux  propri\'et\'es suivantes d'un $A$-sch\'ema $\X$ de type (R).

\bigskip

{\it (i) La fibre sp\'eciale  $Y/F$
    contient une  composante g\'eo\-m\'e\-tri\-quement int\`e\-gre de multiplicit\'e 1.

(ii) La fibre sp\'eciale $Y/F$
    contient une  composante g\'eo\-m\'e\-tri\-quement int\`e\-gre de multiplicit\'e 1
    qui admet un $F$-morphisme depuis une $F$-vari\'et\'e s\'epara\-ble\-ment rationnellement connexe.

 (iii) La fibre sp\'eciale $Y/F$
    contient une  composante g\'eo\-m\'e\-tri\-quement int\`egre de multiplicit\'e 1
    qui est  une $F$-vari\'et\'e s\'epara\-ble\-ment rationnellement connexe.}
    
  On laisse ici au lecteur  le soin de v\'erifier que la propri\'et\'e (ii) satisfait la
  m\^eme propri\'et\'e d'invariance $K$-birationnelle que la propri\'et\'e (i)
  (cf. \S \ref{compgeom1}).

\begin{thm}  (CT-Kunyavski\v{\i} 2006) \cite{CTK}
\label{CTK}
 Soit $A$ un anneau de valuation discr\`ete, de corps des fractions $K$, de corps r\'esiduel $F$
 de caract\'eristique z\'ero. Soit $\X$ un $A$-sch\'ema r\'egulier propre int\`egre de fibre g\'en\'erique $\X_{K}$ une compactification lisse d'un espace  homog\`ene principal d'un groupe 
 semi-simple simplement connexe, \`a fibre sp\'eciale un diviseur \`a croisements normaux stricts.  Il existe alors une composante de la fibre sp\'eciale qui est
 g\'eo\-m\'e\-tri\-quement int\`egre et de multiplicit\'e 1, 
 et qui de plus admet un $F$-morphisme
 depuis une $F$-vari\'et\'e rationnellement connexe.
   \end{thm}

\begin{proof}   Comme rappel\'e au paragraphe \ref{rappels},
on peut suivant Ducros \cite{Ducros} plonger $F$ dans un corps $L$ satisfaisant :

(i) Le corps $F$ est alg\'ebriquement ferm\'e dans $L$.

(ii) Le corps $L$ est un corps de dimension cohomologique 1.

(iii) Le corps $L$ est limite inductive   
de corps de fonctions de $F$-vari\'et\'es admettant des fibrations successives
(par applications rationnelles)
en vari\'et\'es qui sont 
des restrictions \`a la Weil de vari\'et\'es de Severi-Brauer.
On voit ais\'ement que de telles vari\'et\'es sont  birationnelles \`a des vari\'et\'es rationnellement connexes (on n'a pas ici besoin  d'invoquer le th\'eor\`eme \ref{ratconsurratcon}).

D'apr\`es Bruhat et Tits, tout espace  homog\`ene principal sous un groupe semi-simple
simplement connexe sur le corps local $L((t))$, dont  le corps r\'esiduel est parfait
et de dimension cohomologique 1, est trivial, i.e. poss\`ede un point  $L((t))$-rationnel.
 Je renvoie ici le lecteur \`a \cite{CTK}  pour l'alg\`ebre commutative utilis\'ee pour terminer la d\'emonstration.
 \end{proof}

\begin{rmq}
L'assertion sur l'existence d'un 
$F$-morphisme
 depuis une $F$-vari\'et\'e rationnellement connexe  ne
 figurait pas  dans \cite{CTK}.  \end{rmq}

\begin{thm}  \cite{CTnoteC2}
\label{CTnoteC2} 
 Soit $A$ un anneau de valuation discr\`ete de corps des fractions $K$, de corps r\'esiduel $F$
 de caract\'eristique z\'ero.  Soit  \  $\Phi \in A[x_{0},\dots,x_{n}]$ une forme homog\`ene de degr\'e $d$ en $n+1>  d^2$ variables.
   Supposons que  l'hypersurface $X/K$ d\'efinie par $\Phi=0$ dans $\bP^n_{K}$ est lisse. 
   Soit $\X/A$ un  mod\`ele r\'egulier  de cette hypersurface, propre et plat sur $A$, \`a fibre sp\'eciale 
    \`a croisements normaux stricts.
    Il existe alors une composante de la fibre sp\'eciale qui est
 g\'eo\-m\'e\-tri\-quement int\`egre et de multiplicit\'e 1, 
 et qui de plus admet un $F$-morphisme
 depuis une $F$-vari\'et\'e rationnellement connexe.
   \end{thm}

\begin{proof}
 Pour \'etablir le r\'esultat on peut supposer $A=F[[t]]$.
 D'apr\`es le th\'eor\`eme \ref{grandcorpsHogadiXu}, on peut plonger
$F$ dans un corps $L$ qui est union de corps de fonctions
de $F$-vari\'et\'es rationnellement connexes, et qui est un corps $C_{1}$.
On remplace $F[[t]]$ par $L[[t]]$ et on utilise le fait que $L((t))$ est un corps $C_{2}$
 puisque $L$ est un corps $C_{1}$ (th\'eor\`eme \ref{Greenberg}). 
 On a donc $X(L((t))\neq \emptyset$ et donc $\X(L[[t]]) \neq \emptyset$.
 De ceci on d\'eduit que la fibre sp\'eciale contient une composante g\'eom\'etriquement int\`egre
 de multiplicit\'e 1. 
On termine alors la d\'emonstration comme dans le th\'eor\`eme  \ref{CTK} ci-dessus.
 \end{proof}

 \begin{rmq}
L'assertion sur l'existence d'un 
$F$-morphisme
 depuis une $F$-vari\'et\'e rationnellement connexe  
 ne  figurait pas  dans \cite{CTnoteC2}. C'est l'utilisation du
 th\'eor\`eme \ref{HogadiXu}  (Hogadi et Xu) au lieu du th\'eor\`eme \ref{KollarAx}  (Koll\'ar)
 qui permet ici de l'obtenir.
 \end{rmq}

\begin{rmq}
Soit $A$ l'anneau des entiers d'un corps $p$-adique, $\F$ son corps r\'esiduel. 
Soit $ \Phi , n, d$ et $\X/A$ comme dans  l'\'enonc\'e du th\'eor\`eme \ref{CTnoteC2},
 en particulier on suppose donn\'e 
 un mod\`ele \`a fibre sp\'eciale \`a croisements normaux stricts.
Supposons que le th\'eor\`eme \ref{CTnoteC2}  vaille encore dans ce
cas d'in\'egale caract\'eristique.
 
  On dispose alors  d'une composante 
 de $Y$
qui est de multiplicit\'e 1 et est
g\'eom\'etriquement int\`egre  sur le corps fini $\F$. Par les estimations de Lang-Weil,
 il existe un z\'ero-cycle de degr\'e 1 (par rapport au corps $\F$) de support dans le lieu lisse
 de  cette composante et 
  non situ\'e sur les autres composantes.
 Par le lemme de Hensel,
on peut relever ce z\'ero-cycle sur $K$  et l'on obtient que $X/K$ poss\`ede un z\'ero-cycle de
degr\'e 1.  Pour $K$ un corps $p$-adique, l'existence d'un z\'ero-cycle de degr\'e 1 sur toute  hypersurface lisse
 de degr\'e $d$ dans $\bP^n_{K}$, avec $n \geq d^2$,  
est une conjecture de Kato et Kuzumaki  \cite{KK}, \'etablie par ces auteurs  
lorsque $d$ est un nombre premier.

 On dispose plus pr\'ecis\'ement d'une $\F$-application rationnelle
d'une $\F$-vari\'et\'e s\'eparablement rationnellement connexe   sur un corps fini
vers une composante lisse de multiplicit\'e 1 de la fibre sp\'eciale.
Le th\'eor\`eme \ref{Esnault1} (Esnault) assure l'existence d'un $\F$-point sur toute
$\F$-vari\'et\'e s\'eparablement rationnellement connexe, donc par le lemme
\ref{Nishi} sur la composante lisse.
 Mais tout tel $\F$-point peut se trouver aussi sur une autre
composante, donc ne pas \^etre lisse sur $Y$, ce qui emp\^eche de le relever
en un $K$-point de $X$.
C'est heureux. Sinon
 (modulo l'existence de bons mod\`eles) $\Q_{p}$ serait un corps $C_{2}$
 (ex-conjecture d'E. Artin).
Mais les exemples fameux de Terjanian et de ses successeurs montrent
que $\Q_{p}$ n'est pas $C_{2}$. 

Il vaudrait d'ailleurs la peine
de regarder les nombreux contre-exemples  \`a la conjecture d'Artin qui ont \'et\'e construits et de v\'erifier
qu'il existe  toujours dans ces cas un z\'ero-cycle de degr\'e 1. Il en est ainsi
pour l'exemple initial de Terjanian sur $\Q_{2}$.
\end{rmq}

\subsection{Vari\'et\'es rationnellement simplement connexes}

Les vari\'et\'es rationnellement connexes sont un analogue alg\'ebrique des espaces topologiques connexes par arcs. B. Mazur a demand\'e s'il y a un analogue en g\'eom\'etrie alg\'ebrique  des espaces simplement connexes. En topologie, on demande que l'espace des lacets point\'es soit connexe par arcs. A la suite d'une suggestion de Mazur,   de Jong et Starr \cite{dJS2}
 proposent les d\'efinitions suivantes. Dans l'\'etat actuel des recherches, il faut consid\'erer ces d\'efinitions
comme provisoires.

Soit $X$ une vari\'et\'e projective et lisse sur $\C$, \'equip\'ee d'un fibr\'e ample $H$.
Soit ${\overline M}_{0,2}(X,e)$ l'espace de Kontsevich param\'etrisant les donn\'ees suivantes : une courbe $C$ propre, r\'eduite, connexe, \`a croisements normaux, de genre arithm\'etique 0,
un couple ordonn\'e  $(p,q)$ de points lisses de $C$,
un morphisme $h : C \to X$ de cycle image 
 de degr\'e $e$,
tels que de plus la situation n'ait qu'un nombre fini d'automorphismes.

On dispose alors d'un morphisme d'\'evaluation
$$ {\overline M}_{0,2}(X,e)  \to X \times X.$$
La fibre g\'en\'erale de ce morphisme est un analogue de l'espace des chemins
\`a points base en topologie.

La vari\'et\'e (projective et lisse) $X$ est dite {\it rationnellement  simplement connexe}
 si pour $e \geq 1$ suffisamment grand il existe une composante $M$
de ${\overline M}_{0,2}(X,e)$ dominant $X \times X$
telle que la fibre g\'en\'erique de $M \to X \times X$ soit 
une vari\'et\'e rationnellement connexe.

\medskip

 De Jong et Starr  \cite{dJS2} 
  consid\`erent  aussi
l'espace $$ {\overline M}_{0,m}(X,e)$$ o\`u cette fois-ci l'on fixe
$m \geq 2$ points lisses ordonn\'es sur la courbe de genre arithm\'etique  z\'ero, et l'\'evaluation
$$ {\overline M}_{0,m}(X,e)  \to X^m.$$

Ils appellent $X$ {\it fortement rationnellement simplement connexe}
si pour tout $m\geq 2$ et tout entier $e $ suffisamment grand
(fonction de $m$) il existe une composante $M$
de ${\overline M}_{0,m}(X,e)$ dominant $X^m$
telle que la fibre g\'en\'erique de $M \to X^m$ soit 
une vari\'et\'e rationnellement connexe.

\bigskip

De Jong et Starr (travaux en cours) ont obtenu une s\'erie de r\'esultats sur les intersections compl\`etes
lisses dans l'espace projectif. Pour simplifier, je cite leurs r\'esultats pour les hypersurfaces.

\begin{thm} (de Jong-Starr) \cite{dJS2}
\label{hypsimpsimp}
  Une hypersurface lisse de degr\'e $d \geq 2$ dans $\bP^n_{\C}$ avec $$n \geq d^2-1$$
est rationnellement simplement connexe, \`a l'exception des quadriques dans $\bP^3_{\C}$.
\end{thm}

\begin{thm} (de Jong-Starr) \cite{dJS2}
\label{hypfortsimp}
  Une hypersurface lisse de degr\'e $d \geq 2$ dans $\bP^n_{\C}$ avec 
  $$n \geq 2d^2-d-1$$
est fortement rationnellement simplement connexe.
\end{thm}

Dans la d\'efinition ci-dessus on peut prendre $e \geq 4m-6$.

\begin{thm} (de Jong-Starr) \cite{dJS2} Pour  $n \geq d^2$, 
il existe un ouvert de Zariski non vide de l'espace des hypersurfaces de degr\'e $d$
dans  $\bP^n_{\C}$ tel que toute hypersurface param\'etr\'ee par un point de cet espace 
est fortement  rationnellement  simplement connexe.
\end{thm}

 La suggestion suivante est une version locale d'une suggestion globale de de Jong
 (\ref{sugdj} ci-apr\`es).

 \begin{sugn} Soit $A$ un anneau de valuation discr\`ete de corps des fractions $K \subset \C$
et soit $F$ son  corps r\'esiduel, suppos\'e de caract\'eristique z\'ero. Soit $\X$
un $A$-sch\'ema de type (R), $X/K$ sa fibre g\'en\'erique, $Y/F$ sa fibre sp\'eciale. Si les conditions suivantes sont satisfaites  :

(i) la $\C$-vari\'et\'e $X\times_{K}\C$ est fortement rationnellement simplement con\-nexe,

(ii) l'obstruction \'el\'ementaire pour $X/K$ s'annule,
 
\noindent alors la fibre sp\'eciale $Y/F$  contient une  composante g\'eo\-m\'e\-tri\-quement int\`e\-gre de multiplicit\'e 1 qui admet un $F$-morphisme depuis une $F$-vari\'et\'e 
 rationnellement connexe.

\end{sugn}

\begin{rmq} 
Dans (i),  l'exemple \ref{exdj}  et le th\'eor\`eme \ref{hypsimpsimp} justifient la restriction aux vari\'et\'es
fortement rationnellement simplement con\-nexes, plut\^ot qu'aux vari\'et\'es rationnellement simplement con\-nexes.
L'exemple \ref{extao} justifie la condition (ii). 
\end{rmq}

\subsection{Existence d'un point rationnel sur un corps de fonctions de deux variables}
\label{pointrat2}

Sur $K$ un corps de fonctions d'une variable sur $\C$, le th\'eor\`eme de Graber-Harris-Starr 
dit que les $K$-vari\'et\'es rationnellement connexes ont automatiquement un point $K$-rationnel (et le th\'eor\`eme de Graber-Harris-Starr-Mazur dit que ce sont essentiellement les seules).

On peut se demander s'il existe une classe de vari\'et\'es qui ont la propri\'et\'e que lorsqu'elles sont d\'efinies sur un corps $K$ de  fonctions de deux variables sur $\C$, elles ont automatiquement un $K$-point.

Voici deux familles de vari\'et\'es pour lesquelles ceci est connu.

Le th\'eor\`eme de Tsen-Lang implique que  toute hypersurface de degr\'e $d$ dans $\bP^n_{K}$
 avec $n \geq d^2$ poss\`ede un $K$-point. 
 
Soit $G$ un $K$-groupe semi-simple simplement connexe, $E$ un espace homog\`ene principal 
de $G$ et $X$ une $K$-compactification lisse de $E$. 
La conjecture II de Serre pour le corps $K$ affirme que   $E$ et donc aussi $X$
ont un $K$-point. Ceci est connu lorsque $G$ n'a pas de facteur de type $E_{8}$
(Merkur'ev-Suslin, Suslin, Bayer-Parimala, P. Gille).
Pour avoir l'\'enonc\'e dans tous les cas il reste \`a traiter le cas $E_{8}$ d\'eploy\'e.
La r\'esolution de ce dernier cas a \'et\'e  r\'ecemment annonc\'ee par de Jong et Starr,
leur d\'emonstration utilise les techniques 
de vari\'et\'es rationnellement simplement connexes.

\begin{sugn} (de Jong)
\label{sugdj}
Soit $K=\C(S)$ le corps de fonctions d'une surface sur le corps des complexes.
Soit $X$ une $K$-vari\'et\'e 
 fortement rationnellement simplement connexe.
Supposons   l'application de restriction
$\Br K \to \Br K(X)$   injective.
Alors $X$ poss\`ede un point $K$-rationnel.
\end{sugn}

\begin{rmq} L'exemple 
\ref{extao}
montre la n\'ecessit\'e de la condition non g\'eom\'etrique portant sur
le groupe de Brauer.
Des conditions suppl\'e\-mentaires
 de m\^eme nature  
pourraient \^etre n\'ecessaires.
Par exemple on peut demander que pour toute extension finie (ou non)
de corps $L/K$ l'application $\Br L \to \Br L(X)$ soit injective.
 De fa\c con encore plus g\'en\'erale, on peut demander 
que pour toute extension finie (ou non)
de corps $L/K$
 il n'y ait pas  d'obstruction \'el\'ementaire
\`a l'existence d'un $L$-point sur $X\times_{K}L$ (voir  \cite{BoCTSk}).

Pour une $K$-vari\'et\'e $X$  intersection compl\`ete lisse de dimension au moins 3 dans un espace projectif $\bP^n_{K}$, 
l'obstruction \'el\'ementaire s'annule.
 Il en est de m\^eme pour une $K$-vari\'et\'e projective et lisse
g\'eom\'etriquement connexe qui contient un ouvert $U$ qui est un espace homog\`ene principal  
d'un groupe semi-simple simplement connexe. Pour ces r\'esultats, voir  \cite{BoCTSk}.
\end{rmq}

 \begin{rmq}
Dans \cite{BoCTSk} on s'int\'eresse aux compactifications lisses d'espaces homog\`enes de groupes lin\'eaires connexes sur $K=\C(S)$ un corps de fonctions de deux variables, lorsque les stabilisateurs g\'eom\'e\-triques sont connexes
(et qu'il n'y a pas de facteur $E_{8}$). On montre que dans ce cas l'obstruction \'el\'ementaire
\`a l'existence d'un point rationnel est la seule obstruction. 
 
Pour une hypersurface cubique lisse de dimension au moins 3 sur un corps $K$,
l'obstruction \'el\'ementaire s'annule. Sur l'exemple \ref{exdj}
on voit donc que l'obstruction \'el\'ementaire est loin de contr\^oler l'existence
d'un point rationnel pour les vari\'et\'es rationnellement connexes sur un corps de fonctions de deux variables
sur les complexes.
\end{rmq}

\begin{rmq}
De Jong et Starr ont un travail en pr\'eparation sur les vari\'et\'es rationnellement simplement connexes o\`u ils montrent que certains espaces homog\`enes projectifs sur un corps de   fonctions de deux variables sur $\C$  ont automatiquement un point rationnel. Cela leur permet de donner une nouvelle d\'emonstration (la troisi\`eme !) du th\'eor\`eme de de Jong \cite{dJindexexp} qu'indice et exposant co\"{\i}ncident  pour les alg\`ebres simples centrales sur un tel corps.
\end{rmq}

\subsection{Approximation faible en toutes les places d'un corps de fonctions d'une variable}
\label{apfaiblepartout}
 
 Rappelons que c'est une question ouverte (\ref{appfaibleratcon}) de savoir si toute vari\'et\'e rationnellement connexe
sur un corps de fonctions d'une variable satisfait l'approximation faible en toute place.

\begin{thm} (Hassett-Tschinkel)\cite{HT3}
 Soit $K$ un corps de fonctions d'une variable sur un corps alg\'ebriquement clos de caract\'eristique z\'ero.  Il existe une fonction $\varphi : \N \to \N$ satisfaisant la propri\'et\'e suivante.
  Pour toute hypersurface lisse de degr\'e $d$ dans $\bP^n$ avec $n \geq \varphi(d)$, l'approximation faible vaut en tout ensemble fini de places de $K$.
  \end{thm}
  
  Pour $d=3$, $\varphi(3)=6$ convient.

  Un travail en cours  sur les vari\'et\'es  rationnellement  simplement connexes (de Jong-Starr  \cite{dJS2}, appendice de Hassett) 
   donne $\varphi(d) \leq  2d^2-d-1$ et, si l'hypersurface est {\og g\'en\'erale \fg}, $\varphi(d)  \leq d^2$.

\begin{thm} (Hassett)
  Soit $K=\C(C)$ le corps des fonctions d'une courbe. Si $X/K$ est une vari\'et\'e
fortement  rationnellement simplement connexe, alors elle satisfait l'approximation
faible par rapport \`a tout ensemble fini de places de $K$.
\end{thm}

\subsection{$R$-\'equivalence et \'equivalence rationnelle}

Dans la recherche de la bonne d\'efinition de vari\'et\'es  {\og sup\'e\-rieurement \fg} ration\-nel\-lement  connexes, on peut aussi penser \`a des conditions de trivialit\'e de $X(k)/R$ et de $A_{0}(X)$ sur les corps {\og de dimension 2 \fg}, comme les corps $p$-adiques, les corps de fonctions de deux variables sur les complexes, les corps de s\'eries formelles it\'er\'ees $\C((a))((b))$.

\subsubsection{Groupes semi-simples simplement connexes}

Si $K$ est un corps $p$-adique, ou  si $K$ est un corps de fonctions de deux variables sur
les complexes,  ou si $K=\C((a))((b))$,  et si $G/K$ est un groupe semi-simple simplement connexe sans facteur de type $E_{8}$, on sait \'etablir $G(K)/R=1$ et $X(K)/R=1$ (voir \cite{CTGillePa}). Pour $X$ une compactification lisse d'un tel $G$, ceci implique
$A_{0}(X)=0$.

\subsubsection{Hypersurfaces cubiques lisses}

\begin{prop} (Madore) \cite{MadoreJNT}
 Soit $K$ un corps $p$-adique ou un corps $C_{2}$. Soit $X \subset \bP^n_{K}$ une
hypersurface cubique lisse. Pour $n \geq 11$, on a \break ${\rm card} \hskip1mm X(K)/R=1$ et $A_{0}(X)=0$.
\end{prop}

Soit $K$ un corps $p$-adique.
Pour $n=3$, on sait donner des exemples avec $X(K)/R$ et $A_{0}(X)$ d'ordre plus grand que 1.
On ignore ce qui se passe pour $4   \leq n \leq 10$.

Par exemple, qu'en est-il pour l'hypersurface cubique d'\'equation   :
$$x^3+y^3+z^3+pu^3+p^2v^3=0$$
dans $\bP^4_{\Q_{p}}$ ?

Supposons $p \equiv 1$ mod $3$, et soit $a \in \Z_{p}^{\times}$ non cube.
Qu'en est-il pour l'hypersurface
$$x^3+y^3+z^3+p(u_{1}^3+au_{2}^3) +p^2(v_{1}^3+av_{2}^3)=0$$
dans $\bP^6_{\Q_{p}}$ ?

Sur le corps $K=\C((a))((b))$, en utilisant la th\'eorie de l'intersection sur un mod\`ele
au-dessus de $\C((a))[[b]]$,  Madore \cite{MadoreJNT} a  montr\'e que pour l'hypersurface cubique
lisse $X \subset \bP^4_{K}$  
  d'\'equation
$$x^3+y^3+az^3+bu^3+abv^3=0,$$
on a $A_{0}(X) \neq 0.$

\subsubsection{Intersections  lisses de deux quadriques}

Soit $K$ un corps $p$-adique, et soit $X \subset \bP^n_{K}$, avec $n \geq 4$, une intersection compl\`ete  lisse de deux quadriques poss\'edant un $K$-point.  
 Si $n \geq 7$, alors 
${\rm card} \hskip1mm X(K)/R = 1$ (\cite{CTSaSwD}) et donc $A_{0}(X)=0$ (ceci vaut aussi pour un corps $C_{2}$).
L'ensemble  fini $X(K)/R$ peut \^etre non trivial pour $n=4$. Les cas $n=5$ et $n=6$  sont ouverts.
 Le groupe $A_0(X)$ est   nul si $n=6$ et  $k$ est non dyadique (\cite{ParimalaSuresh}).
Le groupe fini $A_{0}(X)$ peut \^etre  non trivial pour $n=4$.
Les cas $n=5$ et $n=6$ ($k$ dyadique) sont ouverts.

Soient $a_{i}, i=1, \dots, 3$, $b_{i}, i=1, \dots, 3$ dans $\Z_{p}$
satisfaisant $a_{i}\neq b_{i}$ et $a_{i}b_{j}-a_{j}b_{i} \in \Z_{p}^{\times}$
pour $i \neq j$.
Soit $X \subset \bP^5_{\Q_{p}}$ l'intersection compl\`ete lisse de deux quadriques
donn\'ee  par le syst\`eme
$$ \sum_{i=0}^3a_{i}X_{i}^2+pX_{4}^2=0, \hskip2mm  \sum_{i=0}^3b_{i}X_{i}^2+pX_{5}^2=0.$$
Que valent $X(\Q_{p})/R$ et $A_{0}(X)$ ?

\subsubsection{Fibr\'es en quadriques sur la droite projective}

Soit $K$ un corps $p$-adique, et soit $X $ une $K$-vari\'et\'e g\'eom\'etriquement int\`egre, projective et lisse sur $K$, fibr\'ee en quadriques de dimension $d\geq 1$ sur la droite projective $\bP^1_{K}$.
Si $p\neq 2$ et $d\geq 3$, alors $A_{0}(X)=0$ (\cite{ParimalaSuresh}). Dans le cas $d=2$,
Parimala et Suresh \cite{ParimalaSuresh}
  ont un exemple   int\'eressant avec $A_{0}(X) \neq 0$.  Dans cet exemple, un \'el\'ement non nul de $A_{0}(X) $  est d\'etect\'e par la mauvaise r\'eduction de $X$ de fa\c con subtile,   le groupe de Brauer de $X$
ne permet pas   de d\'etecter cet \'el\'ement.

\section{Surjectivit\'e arithm\'etique et surjectivit\'e \\ g\'eom\'etrique}

Mis \`a part bien s\^ur le th\'eor\`eme d'Ax et Kochen, les \'enonc\'es de ce paragraphe sont \'etablis dans des notes non publi\'ees de l'auteur. Le lecteur ne devrait pas avoir de difficult\'e \`a reconstituer les d\'emonstrations.

\subsection{Morphismes d\'efinis sur un corps de nombres et applications induites sur les points 
locaux}

On s'int\'eresse dans la suite \`a la situation suivante.
 
 \medskip
 
 (*) On est sur un corps de nombres $k$, $X$ et $Y$ sont deux $k$-vari\'et\'es lisses
 g\'eo\-m\'e\-tri\-quement int\`egres, la $k$-vari\'et\'e $Y$ est projective,
 on a un $k$-morphisme  projectif  $ f : X \to Y$ de fibre g\'en\'erique g\'eom\'etriquement int\`egre.
 On note $U \subset X$ l'ouvert de lissit\'e du morphisme $f$.

 \medskip

 On demande quels sont les liens entre la g\'eom\'etrie du morphisme $f$ et 
 les propri\'et\'es de surjectivit\'e des applications induites  $X(k_{v}) \to Y(k_{v})$ pour presque toute place  $v$, ou d\'ej\`a pour une infinit\'e de places $v$ du corps de nombres $k$.

Les th\'eor\`emes \ref{surjinftv}  et \ref{surjpptv} ci-apr\`es jouent un   r\^ole central dans l'\'etude du
principe de Hasse pour les vari\'et\'es alg\'ebriques sur un corps de nombres.

\begin{thm} 
\label{surjinftv}
Sous les hypoth\`eses (*), si $Y$
 est une {\it courbe} et
 si l'application induite $U \to Y$ est surjective  (ce qui \'equivaut \`a :
$f : X \to Y$ est localement scind\'e pour la topologie \'etale sur $Y$),
alors   il existe une infinit\'e de places
$v$ de
$k$ pour lesquelles l'application induite $X(k_{v}) \to Y(k_{v})$
est surjective.
\end{thm} 

\begin{proof}    Pour chaque point ferm\'e $P$ de $Y$
 \`a fibre $X_{P}=f^{-1}(P)$ non lisse
 on choisit une composante $Z_{P}$ de multiplicit\'e 1   de   $X_{P}$.
 L'existence d'une telle composante est garantie par l'hypoth\`ese
 que la fibration est localement scind\'ee pour la topologie \'etale sur $Y$.
 Soit $k_{P}$ le corps r\'esiduel de $Y$ en $P$.
 Soit $K_{P}$ la cl\^oture int\'egrale de $k_{P}$ dans le corps des fonctions
 de $Z_{P}$. Soit $K/k$ une extension finie galoisienne de $k$ dans laquelle se plongent
 toutes les extensions $K_{P}/k$.  La fibration $f_{K} : X_{K} \to Y_{K}$ satisfait alors
 les hypoth\`eses du th\'eor\`eme \ref{surjpptv} ci-apr\`es. En combinant ce th\'eor\`eme
 et le  th\'eor\`eme de  Tchebotarev, qui garantit l'existence d'une infinit\'e de places
 $v$ de $k$ d\'ecompos\'ees dans $K$, on conclut. 
Cette d\'emonstration montre que l'ensemble infini de places cherch\'e contient
un ensemble de places de $k$ de densit\'e positive.
 \end{proof}

 \begin{rmqs}
 
 (1)    
  L'hypoth\`ese que l'application $f : X \to Y$ est localement scind\'ee pour la topologie
\'etale sur $Y$ est en particulier satisfaite si apr\`es extension finie convenable de $k$
la fibration $f $ admet une section.  
D'apr\`es le th\'eor\`eme \ref{GHSJ} (Graber, Harris et Starr), c'est   le cas
 si la fibre g\'en\'erique est une vari\'et\'e rationnellement connexe.

(2)
 Le th\'eor\`eme   ne s'\'etend pas \`a $Y$ de dimension sup\'erieure,
 comme l'on voit en consid\'erant 
 une fibration en coniques sur $\bP^2_{\Q}$ dont le lieu de
ramification est une courbe  $C$ lisse et dont  
le rev\^etement double   $D \to C$ associ\'e est donn\'e par une courbe $D/ \Q$ g\'eo\-m\'e\-tri\-quement int\`egre.
 On peut par exemple prendre pour $C \subset \bP^2_{\Q}$ une  courbe elliptique d'\'equation affine $v^2=u(u-a)(u-b)$ et
une famille de coniques d'\'equation g\'en\'erique
$$X^2- uY^2-(v^2-u(u-a)(u-b))T^2=0.$$
\end{rmqs}
En utilisant le th\'eor\`eme de Lang-Weil, 
 on \'etablit le th\'eor\`eme suivant.

\begin{thm} 
\label{surjpptv}
Sous les hypoth\`eses (*),
s'il existe un ouvert $V \subset X$ tel que le morphisme induit $V \to Y$
soit lisse, surjectif,  \`a  fibres g\'eo\-m\'e\-tri\-quement int\`egres, alors pour presque
toute place $v$ de $k$, l'application induite
$X(k_v) \to Y(k_v)$ est surjective. 
 \end{thm}

\begin{rmq}
\label{pascodim1}
Il ne suffit pas d'avoir la propri\'et\'e
en codimension 1 sur $Y$, comme le montre l'exemple suivant.
Prendre $a \in \Q$ non carr\'e et $X \subset \bP^3 \times_{\Q} \bP^2$ donn\'ee
par
$$uX_{0}^2 -avX_{1}^2+wX_{2}^2-a(u+v+w)X_{3}^2=0.$$
Pour une infinit\'e de $p$, la fl\`eche $X(\Q_{p}) \to \bP^2(\Q_{p})$ n'est pas surjective :
pour $a \in \Z_{p}$, $a$ non carr\'e dans $\Z_{p}$,  et $M$ un point
$(p^{2n+1}\alpha,p^{2m+1}\beta,1)$ avec $n,m \geq 0$, $\alpha$ et $\beta$ in $\Z_{p}^*$
et $\alpha.\beta$ carr\'e dans $\Z_{p}$, la fibre en $M$ n'a pas de $\Q_{p}$-point.
\end{rmq}
  
\bigskip

Le c\'el\`ebre th\'eor\`eme d'Ax et Kochen  \cite{AK} peut se  formuler  de la fa\c con suivante.

\begin{thm} (Ax et Kochen) 
Fixons des entiers $d$ et $n  \geq d^2$. Soit $N+1$ la dimension de l'espace des formes homog\`enes de
degr\'e $d$ en $n+1$ variables. Soit $F(x_{0},\dots,x_{N}; y_{0}, \dots, y_{n})$ la forme universelle
de degr\'e $d$ en $n+1$-variables.  Soit $Z \subset \bP^N \times_{\Q} \bP^n$ le ferm\'e d\'efini par l'annulation de cette forme. Soit $\pi : Z \to \bP^N$ la projection sur le premier facteur. 
Sur tout corps de nombres $k$, pour presque toute place $v$ de $k$,   la projection induite $Z(k_{v}) \to \bP^N(k_{v})$ est surjective.
\end{thm}

\begin{rmq}
En combinant le th\'eor\`eme \ref{surjpptv} et le th\'eor\`eme \ref{CTnoteC2},
on \'etablit un \'enonc\'e du type Ax-Kochen pour la restriction de $Z \to \bP^N$ 
au-dessus d'une droite de $\bP^N$ (passant par un point \`a fibre lisse). 

Si l'on pouvait r\'epondre par l'affirmative \`a la question \ref{bonmodele},
la combinaison du  th\'eor\`eme \ref{CTnoteC2} et du  th\'eor\`eme \ref{surjpptv} donnerait
une nouvelle d\'emons\-tra\-tion du th\'eor\`eme d'Ax et Kochen.

Sans r\'epondre \`a la question  \ref{bonmodele}, Jan Denef a   tout r\'ecemment (juin 2008)
obtenu une nouvelle d\'emonstration du th\'eor\`eme  d'Ax et Kochen, en \'etablissant
 une conjecture g\'en\'erale de  \cite{CTnoteC2}.
  \end{rmq}

On s'int\'eresse aux r\'eciproques des \'enonc\'es ci-dessus.

\medskip

\begin{thm} 
\label{reciprinfv}
 Pla\c cons-nous sous les hypoth\`eses (*), avec $Y$
une courbe.

(i) Si pour une infinit\'e de places $v$  l'application $X(k_{v}) \to Y(k_{v})$ est surjective,
alors pour tout point $P \in Y(k)$ il existe une composante de multiplicit\'e 1
de $f^{-1}(P)$.

(ii) Si pour toute extension finie $K/k$, pour une infinit\'e de places $w$ de $K$,
l'application $X(K_{w}) \to Y(K_{w})$ est surjective, alors l'application induite $U \to Y$ est surjective : le morphisme $X \to Y$
est localement scind\'e pour la topologie \'etale sur $Y$.
\end{thm}

\begin{thm} 
\label{reciprpptv}
 Pla\c cons-nous sous les hypoth\`eses (*), avec $Y$
une courbe.
 Supposons que  pour presque toute place $v$ de $k$ l'application $X(k_{v}) \to Y(k_{v})$ est surjective.
 
Alors :

(a) L'application induite $U \to Y$ est surjective : le morphisme $X \to Y$
est localement scind\'e pour la topologie \'etale sur $Y$. 

(b) Toute fibre connexe
de $U \to Y$ est g\'eo\-m\'e\-tri\-quement connexe. 

(c) Si $P$ est un point ferm\'e de $Y$,
de corps r\'esiduel $\kappa$
et la fibre $f^{-1}(P)$ s'\'ecrit $\sum_{i}e_{i}D_{i}$ avec chaque $D_{i}$ diviseur int\`egre,
de corps des fonctions $\kappa_{i}$, la fl\`eche
$$H^1(\kappa,\Q/\Z) \to \oplus_{i} H^1(\kappa_{i},\Q/\Z)$$ 
obtenue par somme des applications $e_{i}. {\rm Res}_{\kappa_{i}/\kappa}$ est injective.
\end{thm}

\begin{rmq}

Le th\'eor\`eme \ref{reciprpptv}  est le meilleur possible, comme le montre l'exemple
suivant.
Soit $k$ un corps de nombres, $a, b  \in k^*$ avec $a, b , ab \notin k^{*2}$.
Soit $f : X \to \bP^1_{k}$ un mod\`ele projectif de la situation affine suivante :
$$(x^2-ay^2)(u^2-bv^2)(z^2-abw^2)=t,$$
la fl\`eche de projection sur $\A^1_{k}$ \'etant donn\'ee par la coordonn\'ee $t$.
Alors
\medskip

a) La fibre de $f$ en $0$ ne contient aucune composante g\'eo\-m\'e\-tri\-quement int\`egre
de multiplicit\'e 1.

b) Pour toute extension finie $K/k$, pour presque toute place $w$ de $K$, l'application $X(K_{w}) \to Y(K_{w})$ est surjective.

\medskip

Tout le probl\`eme est qu'un polyn\^ome en une variable sur un corps de nombres peut
avoir une solution partout localement sans en avoir sur le corps de nombres, d\`es qu'il
est r\'eductible.
 Ainsi on ne peut pas partir du th\'eor\`eme d'Ax et Kochen pour en d\'eduire le 
th\'eor\`eme \ref{CTnoteC2}
 sur la r\'eduction des formes lisses de degr\'e $d$ en $n> d^2$ variables.

\end{rmq}

 \subsection{Quelques autres questions}

Diverses questions connexes ont \'et\'e discut\'ees dans la litt\'erature.

\medskip

Soient $k$ un corps de nombres et
$f : X \to Y$ un $k$-morphisme propre de $k$-vari\'et\'es lisses g\'eo\-m\'e\-tri\-quement int\`egres,
  $Y$ \'etant une courbe.

\medskip

Question 1. Si sur toute extension finie $K/k$ l'application $X(K) \to Y(K)$ est surjective (\`a un nombre fini
de points pr\`es),
le morphisme admet-il une section ?

\medskip

Question  2. Si sur tout compl\'et\'e $k_{v}$   le $k_{v}$-morphisme  $X_{k_{v}} \to Y_{k_{v}}$ a une section,
le morphisme $f$  admet-il une section ?

\medskip

En dimension relative 1,  pour les courbes relatives de genre z\'ero,
la r\'eponse \`a ces deux questions est  oui  pour $Y=\bP^1$ (Schinzel, Salberger, Serre).
Ceci utilise l'injection 
$\br k(t) \hookrightarrow \prod_{v} \br k(t) $
qui s'\'etablit en consid\'erant la suite exacte de localisation pour le groupe de Brauer
sur la droite projective.
La r\'eponse est non pour  $Y$ une courbe de genre 1 : on utilise une courbe
elliptique $Y$ avec un \'el\'ement de 2-torsion dans son groupe de Tate-Shafarevich 
repr\'esent\'e par une alg\`ebre de quaternions sur le corps de fonctions $k(Y)$ (voir \cite{ParimalaSujatha}).

\medskip

En dimension relative 1, pour les courbes relatives de genre 1 et $Y=\bP^1$, la question 1 est ouverte.
La question 2 a une r\'eponse n\'egative (prendre   $X=C\times_{k}\bP^1$ 
avec $C$ une courbe de genre 1 qui est un contre-exemple au principe de Hasse).

\medskip

Pour les familles de quadriques de dimension relative $d \geq 2$ au-dessus de $Y=\bP^1_{\Q}$
la  r\'eponse aux deux questions ci-dessus est n\'egative, et ce pour tout tel $d$.
 
\medskip

Soit $k$ un corps de nombres totalement imaginaire. 
Pour les familles de quadriques de dimension relative $d \geq 2$ au-dessus de $Y=\bP^1_{k}$,
 la r\'eponse aux deux questions est n\'egative pour $2 \leq d \leq 6$. 
 (Pour $d\geq 7$, on conjecture qu'il y a toujours une section.)
 
 \medskip

  Pour justifier ces r\'eponses n\'egatives,
 partons d'un couple de formes
  quadratiques $f(x_{1}, \dots,x_{n}), g(y_{1}, \dots, y_{m})$ sur le corps
  de nombres $k$ tel que sur tout compl\'et\'e de $k$ l'une des deux formes
   ait un z\'ero non trivial (donc $n$ et $m$ sont au moins \'egaux \`a 2)
   mais que pour chacune de ces formes il existe un compl\'et\'e $k_{v}$
   sur lequel la forme n'a pas de z\'ero non trivial.
   
   Un th\'eor\`eme d'Amer et de Brumer (see the references in \cite{CTSaSwD}) 
   garantit que sur toute extension $F$ de $k$,
   la forme quadratique $f(x_{1}, \dots,x_{n})+tg(y_{1}, \dots, y_{m})$ sur le corps $F(t)$
admet un z\'ero non trivial sur $F(t)$ si et seulement si le syst\`eme
$$f(x_{1}, \dots,x_{n})=0,  \hskip1mmg(y_{1}, \dots, y_{m})=0$$ 
admet un z\'ero non trivial dans $F^{n+m}$.

La forme $f(x_{1}, \dots,x_{n})+tg(y_{1}, \dots, y_{m})$ sur le corps $k(t)$
a alors un z\'ero sur chaque $k_{v}(t)$ mais n'en a pas sur $k(t)$. Ceci donne
les r\'eponses n\'egatives \`a la question 2, et les r\'eponses n\'egatives \`a la
question 1 r\'esultent du principe de Hasse pour les formes quadratiques sur
un corps de nombres.

\vskip2cm 
Jean-Louis Colliot-Th\'el\`ene

Math\'ematiques, B\^atiment 425

Universit\'e  Paris-Sud

F-91405 Orsay FRANCE

jlct \`a math.u-psud.fr


\begin{thebibliography}{foo}

 \bibitem{AraujoKollar} C. Araujo et J. Koll\'ar, Rational curves on varieties, in 
{\it  Higher dimensional varieties and rational points}, Bolyai Society Mathematics Studies {\bf 12},
 Springer, 2003,  13-92.
 
  \bibitem{Ax} J. Ax, The elementary theory of finite fields, Ann. of Math. (2) {\bf 88} (1968) 239-271.
 
\bibitem{AK} J. Ax et S. Kochen, Diophantine problems over local fields, I, Amer. J. Math. {\bf 87} (1965) 605--631.

\bibitem{BoCTSk}  M. Borovoi,  J.-L. Colliot-Th\'el\`ene  et A. N. Skorobogatov, The elementary obstruction and homogeneous spaces, Duke Math. J. {\bf 141} (2008) 321--364.
 

\bibitem{BLR} S. Bosch, W. L\"utkebohmert et M. Raynaud, {\it N\'eron Models}, 
Ergebnisse der Math. und ihrer Grenzg. 3. Folge Band {\bf 21}, Springer-Verlag.

\bibitem{Camp} F. Campana,  Connexit\'e rationnelle des vari\'et\'es de Fano, Annales
scientifiques de l'\'Ecole Normale Sup\'erieure {\bf 25} (1992) 539--545.

 \bibitem{CL} A.  Chambert-Loir, Points rationnels et groupes fondamentaux  : applications de la cohomologie $p$-adique, S\'eminaire Bourbaki, 55e ann\'ee, 2002-2003, no 914.
 
 \bibitem{CTinvmath83} J.-L. Colliot-Th\'el\`ene,  Hilbert's theorem 90 for $K_2$, with application to the Chow groups of rational surfaces,  Invent. math. {\bf 71} (1983) 1-20.

\bibitem{CTinvmath06} J.-L. Colliot-Th\'el\`ene, Un th\'eor\`eme de finitude pour le groupe de Chow des z\'ero-cycles d'un groupe alg\'ebrique lin\'eaire  sur un corps $p$-adique,   Invent. math {\bf  159} (2005) 589-606.

\bibitem{CTnoteC2} J.-L. Colliot-Th\'el\`ene, Fibres sp\'eciales des hypersurfaces de petit degr\'e,
C. R. Acad. Sc. Paris  {\bf 346} (2008) 63-65.

\bibitem{CTGille} J.-L. Colliot-Th\'el\`ene et P. Gille,
Remarques sur l'approximation faible sur un  corps de fonctions d'une variable, in {\it Arithmetic of higher dimensional arithmetic varieties} (ed. B. Poonen et Yu. Tschinkel),  Birkh\"auser, Progress in Mathematics,  2003, 121-133.
 

\bibitem{CTGillePa} J.-L. Colliot-Th\'el\`ene, P. Gille et R. Parimala, Arithmetic of linear algebraic groups over two-dimensional fields,  Duke Math. J.  {\bf 121}, 2004, 285--341. 
 


 \bibitem{CTK} J.-L. Colliot-Th\'el\`ene et B. Kunyavski\u{\i}, Groupe de Picard et groupe de Brauer des compactifications lisses d'espaces homog\`enes.   J.  Algebraic Geom.  {\bf 15}  (2006) 
 733--752.
 
 
 \bibitem{CTMadore} J.-L. Colliot-Th\'el\`ene et D. Madore, Surfaces de Del Pezzo sans point rationnel sur un 
corps de dimension cohomologique un,  
Journal de l'Institut Math\'ematique de Jussieu {\bf 3} (2004) 1--16.

\bibitem{CTSaito} J.-L. Colliot-Th\'el\`ene et S. Saito, Z\'ero-cycles sur les vari\'et\'es $p$-adiques et groupe de Brauer, IMRN (1996), no. 4, 151--160.

 \bibitem{RET} J.-L. Colliot-Th\'el\`ene et J.-J. Sansuc,  
 La $R$-\'equivalence sur les tores, Ann. Sc. E. N. S.  {\bf 10} (1977) 175--229.
 
 
 \bibitem{CTSaSwD} J.-L. Colliot-Th\'el\`ene,  J.-J. Sansuc et Sir Peter Swinnerton-Dyer,
 Intersections of two quadrics and Ch\^{a}telet surfaces, I, J. f\"ur die reine und angew. Math.  (Crelle) {\bf 373} (1987) 37-107; II, ibid. {\bf 374}  (1987) 72-168.
 
 \bibitem{CTSkMan} J.-L. Colliot-Th\'el\`ene et A. N. Skorobogatov,  
$R$-equivalence on conic bundles of degree 4, Duke Math. J. {\bf 54} (1987) 671-677.
 
  \bibitem{CTSko} J.-L. Colliot-Th\'el\`ene et A. N. Skorobogatov,  
  Groupes de Chow des z\'ero-cycles des fibr\'es en quadriques, K-theory {\bf 7} (1993) 477--500.
 

\bibitem{coraytsfasman} D. F. Coray and M. A. Tsfasman, Arithmetic on singular Del Pezzo surfaces.
Proc. London Math. Soc. (3) {\bf 57} (1988), no. 1, 25--87. 

 \bibitem{Debarre} O. Debarre, {\it Higher-dimensional algebraic geometry}, Universitext, Springer, 2001.
 
 
 \bibitem{DJL} J.  Denef,  M. Jarden et  D. J. Lewis, On Ax-fields which are $C_{i}$, 
   Quart. J. Math. Oxford Ser. (2)  {\bf 34}  (1983), no. 133, 21--36.  
   

    
\bibitem{Ducros} A. Ducros, Dimension cohomologique et points rationnels sur les courbes,    J. Algebra  {\bf 203}  (1998)  
349--354. 

\bibitem{Ducros2} A. Ducros, Points rationnels sur la fibre sp\'eciale d'un sch\'ema au-dessus d'un anneau de valuation, Math. Z. {\bf 238}  (2001) 177--185.

 

\bibitem{Esnault1} H. Esnault, Varieties over a finite field with trivial Chow group of $0$-cycles have a rational point, Invent. math. {\bf 151} (1)  (2003) 187--191. 	 


\bibitem{Esnault2} H. Esnault,  Deligne's integrality theorem in unequal characteristic and rational points
over finite fields (with an appendix by H. Esnault and P. Deligne),  Annals of Math. {\bf 164} (2006) 715-730.
 

\bibitem{Esnault3} H. Esnault, Coniveau over $p$-adic fields and points over finite fields,
C. R. Acad. Sc. Paris S\'er. I {\bf 345} (2007) 73-76.

\bibitem{Esnault4} H. Esnault et C. Xu, Congruence for rational points over finite fields
and coniveau over local fields, {\tt arXiv:0704.1273v1}, to appear in Transactions AMS. 
 

\bibitem{Euler} L. Euler, Demonstratio theorematis Fermatiani omnem numerum sive integrum sive fractum
esse summam quatuor pauciorumve quadratorum, N. Comm. Ac. Petrop. 5 (1754/5), 1760, p. 13-58.
 E.242, Opera omnia,  Birkh\"auser 2003, vol. I.2, 338-372.
 
 


\bibitem{FR} N. Fakhruddin et C. S. Rajan, Congruences for rational points on varieties over
finite fields, Math. Annalen {\bf 333} (2005) 797-809.
 

\bibitem{GilleIHES} P. Gille, La R-\'equivalence sur les groupes alg\'ebriques r\'eductifs d\'efinis sur un corps global.   Inst. Hautes \'Etudes Sci. Publ. Math.  {\bf 86}  (1997), 199--235.
 

\bibitem{GilleTAMS} P. Gille, Sp\'ecialisation de la R-\'equivalence pour les groupes r\'eductifs, Trans. Amer. Math. Soc. {\bf 356} (2004),  4465-4474.



\bibitem{GHMS}  T.  Graber, J. Harris, B. Mazur et J. Starr,  Rational connectivity and
sections of families over curves, Ann. sci. \'Ecole Norm. Sup. (4) {\bf 38} (2005) 671-692.

\bibitem{GHS}  T. Graber, J. Harris,  et J. Starr, Families of rationally connected varieties,
J. Amer. Math. Soc. {\bf 16} (2003) 57-67.

 

 \bibitem{Greenberg} M. J. Greenberg, Rational points in Henselian discrete valuation rings, Publications math\'ematiques de l'I.H.\'E.S. {\bf 31} (1966) 59--64.
 
 \bibitem{GrBr} A. Grothendieck, Le groupe de Brauer I, II, III, in 
{\it Dix expos\'es sur la cohomologie des sch\'emas},  North-Holland, Amsterdam,  Masson, Paris (1968).
 
 
 \bibitem{HT1} B. Hassett et Yu. Tschinkel,
 Weak approximation over function fields,   Inventiones Math.,  {\bf 163}, no. 1, 171-190, (2006)
 
 
 \bibitem{HT2} B. Hassett et Yu. Tschinkel,
 Approximation at places of bad reduction for rationally connected varieties,   Pure and Applied Math Quarterly, Bogomolov Festschrift, {\bf 4}, no. 3, 1-24, (2008).

 \bibitem{HT3} B. Hassett et Yu. Tschinkel, Weak approximation for hypersurfaces of low degree,   22 pp., (2006).
 
 
 
 \bibitem{HogadiXu} A. Hogadi et Ch. Xu, Degenerations of rationally connected varieties,
{\tt arXiv:math/ 060666v3} (Avril 2007), to appear in Trans. A. M. S. 
 
 
  \bibitem{dJindexexp} A. J. de Jong, The period-index problem for the Brauer group of an algebraic surface, Duke Math. J. {\bf 123} (2004) 71-94.

\bibitem{dJS} A. J. de Jong et J. Starr, Every rationally connected variety over the function field of a curve has a rational point, Amer. J. Math. {\bf 125} (2003) 567-580.

\bibitem{dJS2}  A. J. de Jong et J. Starr, Low degree complete intersections are rationally simply connected, pr\'epublication.

\bibitem{Kahn} B. Kahn, Zeta functions and motives, Pure Appl. Math. Quarterly, {\bf 5} (1),  \`a  para\^{\i}tre.
   
 \bibitem{KK}   K. Kato et T. Kuzumaki, The dimension of fields and algebraic K-theory, Journal of Number Theory
   {\bf 24} (1986) 229--244.
   
   \bibitem{Kollarlivre} J. Koll\'ar, {\it Rational curves on algebraic  varieties}, Ergebnisse der Mathematik und ihrer Grenzgebiete, 3. Folge, Band 32 Springer-Verlag (1996, r\'e\'edition avec corrections 1999).
   
   \bibitem{Kollarannmath} J. Koll\'ar, Rationally connected varieties over local fields,  Ann. of Math. (2)  150  (1999),  no. 1, 357--367.
   
  
  \bibitem{Kollarjap} J. Koll\'ar, Specialization of zero cycles,  Publ. Res. Inst. Math. Sci.  {\bf 40}  (2004),  no. 3, 689--708.
  
  \bibitem{KollarAx} J. Koll\'ar, A conjecture of Ax and degenerations of Fano varieties,  
  Israel J. of Math. {\bf 162}  (2007) 235--252.
  
     \bibitem{DK} J. Koll\'ar, Looking for rational curves on cubic hypersurfaces, in {\it Higher-Dimensional Geometry over Finite Fields},
Volume 16 NATO Science for Peace and Security Series: Information and Communication Security,
ed. D. Kaledin and Y. Tschinkel (2008). Voir aussi :
     U. Derenthal and J. Koll\'ar, Looking for rational curves on cubic hypersurfaces, matharXivv:0710.5516.
 
  
  \bibitem{KMM} J. Koll\'ar, Y. Miyaoka et S. Mori, Rational connectedness and boundedness of Fano manifolds.  J. Differential Geom.  {\bf 36}  (1992),  no. 3, 765--779.
  
     \bibitem{KollarSzabo} J. Koll\'ar et E. Szab\' o,  Rationally connected varieties over finite fields,  Duke Math. J.  {\bf 120}  (2003),  no. 2, 251--267.


  \bibitem{Lafon}   G. Lafon, Une surface d'Enriques sans point sur $\C((t))$,  C.R. Math. Acad. Sci. Paris {\bf 138} 
    (1) (2004)
 51-54.
  
 
   \bibitem{L}  S. Lang, On quasi algebraic closure,  Ann. of Math. (2)  {\bf 55}  (1952)  373--390.
   
   \bibitem{Madorethese} D. Madore, 
   Sur la sp\'ecialisation de la R-\'equivalence, in 
   {\it Hypersurfaces cubiques, R-\'equivalence et approximation faible}, th\`ese de doctorat, Universit\'e Paris-Sud, 2005.
   
   \bibitem{Madoremanmat} D. Madore, \'Equivalence rationnelle sur les hyperrsurfaces cubiques sur les
   corps $p$-adiques, manuscripta math. {\bf 110} (2003) 171--185.
   
   
   \bibitem{MadoreSMF} D. Madore,  Approximation faible aux places de bonne r\'eduction sur les
   surfaces cubiques sur les corps de fonctions, Bull. Soc. math. France {\bf 134} (4) (2006) 475--485.
   
   \bibitem{MadoreJNT}  D. Madore,  \'Equivalence rationnelle sur les hypersurfaces cubiques de mauvaise r\'eduction,
   J. Number Theory {\bf 128} (2008) 926--944.
   
   \bibitem{Manin} Yu. I. Manin, {\it Cubic forms. Algebra, geometry, arithmetic}. Translated from the Russian by M. Hazewinkel. Second edition. North-Holland Mathematical Library, 4. North-Holland Publishing Co., Amsterdam, 1986.
   
   \bibitem{Merkurjev} A. S. Merkur'ev, R-equivalence on three-dimensional tori and zero-cycles,
   Algebra and Number Theory {\bf 2} (2008) 69--89.
   
\bibitem{ParimalaSujatha}  R. Parimala and R. Sujatha,  
Hasse principle for Witt groups of function fields with special reference to elliptic curves. With an appendix by J.-L. Colliot-Th\'el\`ene.  Duke Math. J.  85  (1996),  no. 3, 555--582.
   
   \bibitem{ParimalaSuresh} R. Parimala et V. Suresh, Zero-cycles on quadric fibrations: Finiteness theorems and the cycle map, Invent. math. {\bf 122} (1995) 83-117.
   
   \bibitem{Pf} A. Pfister, 
 {\it  Quadratic forms with applications to algebraic geometry and topology}, London Mathematical Society Lecture Note Series  {\bf 217}, Cambridge University Press, Cambridge, 1995.  
  
  \bibitem{SS07} S. Saito et K. Sato, A finiteness theorem for zero-cycles over $p$-adic fields,
  \`a para\^{\i}tre dans Annals of Math.

 
\bibitem{SerreCG}
J-P.~Serre, {\it Cohomologie galoisienne}, Cinqui\`eme \'edition,
r\'evis\'ee
et compl\'et\'ee, Springer Lecture Notes in Mathematics {\bf 5} (1994)
 

\bibitem{Starr} 
J. Starr, Degenerations of rationally connected varieties and PAC fields, arXiv preprint.

\bibitem{Starrrapport}
J. Starr, Arithmetic over function fields, Clay Mathematics Institute 2006 Summer School on Arithmetic Geometry, pr\'epublication.


\bibitem{Tao} D. Tao, A variety associated to an algebra with involution.  J. Algebra  168  (1994),  no. 2, 479--520.


\bibitem{WittS} O. Wittenberg, La connexit\'e rationnelle en arithm\'etique, notes pour un mini-cours,
Session SMF \'Etats de la Recherche  {\og Vari\'et\'es rationnellement connexes \fg}, Strasbourg, mai 2008.

 
\end{thebibliography}
\end{document}